\RequirePackage{fix-cm}
\documentclass[smallextended]{svjour3}       
\smartqed  
\usepackage{graphicx}
\usepackage{epstopdf}
\usepackage{amssymb}
\usepackage{stfloats}

\usepackage{amsmath}
\usepackage{url}
\usepackage[colorlinks=true,citecolor=blue,linkcolor=blue]{hyperref}
\usepackage{multirow}
\usepackage{booktabs}
\begin{document}

\title{Differential quadrature method for space-fractional diffusion equations on 2D irregular domains}

\titlerunning{Numer Algor}        

\author{X. G. Zhu  \and Z. B. Yuan \and F. Liu \and Y. F. Nie}

\authorrunning{Numer Algor} 

\institute{  X. G. Zhu  \and Z. B. Yuan \and Y. F. Nie \at
               Department of Applied Mathematics, Northwestern  Polytechnical University, Xi'an, Shaanxi 710129, P. R. China \\   
              \email{yfnie@nwpu.edu.cn}       
              \and
               F. Liu \at
               School of Mathematical Sciences, Queensland University of Technology, GPO Box 2434, Brisbane, Qld. 4001, Australia \\
              \email{f.liu@qut.edu.au}          
}

\date{Received: date / Accepted: date}

\maketitle

\begin{abstract}
In mathematical physics, the space-fractional diffusion equations are of particular interest in the studies of physical
phenomena modelled by L\'{e}vy processes, which are sometimes called super-diffusion equations.
In this article, we develop the differential quadrature (DQ) methods for solving the 2D space-fractional diffusion equations on irregular domains.
The methods in presence reduce the original equation into a set of ordinary differential equations (ODEs)
by introducing valid DQ formulations to fractional directional derivatives based on the functional values at scattered nodal points on problem
domain.  The required weighted coefficients are calculated by using radial basis functions (RBFs) as trial functions,
and the resultant ODEs are discretized by the Crank-Nicolson scheme. The main advantages of our methods lie in their flexibility and applicability
to arbitrary domains. A series of illustrated examples are finally provided to support these points.
\keywords{differential quadrature (DQ) \and radial basis functions (RBFs) \and fractional directional derivatives \and space-fractional diffusion equations}
 \subclass{35R11 \and 65D25 \and 65M99}
\end{abstract}

\section{Introduction}\label{intro}
During recent decades, a bulk of attention has been attracted to a special part of partial differential equations (PDEs), i.e., 
the so-called fractional PDEs, especially those abstracted from practical problems. As a new class of mathematical models,
fractional PDEs show good promise in characterizing the physical processes with memory effort and historical dependence, such as the
anomalous dispersion in complex heterogeneous aquifer and fractal geometry \cite{Ref02,Ref03},
thereby making up for the defect of inconsistency with the reality when an integer-order model is applied to describe the similar non-classical phenomena.
The relevant applications include crystal dislocation, hydrodynamics, electrochemistry, plasma turbulence, continuum mechanics,
and so forth \cite{Ref10,Ref08,Ref11,Ref12,Ref09}. There is still a rapidly growing interest on these subjects in this moment. Due to the universal mutuality, however,
their solutions can rarely be represented by elementary functions in closed forms. A few existent analytic solutions are limited
to very simple cases or achieved in series or integral forms under theoretical restrictions;
see \cite{R28,R27,R25,R26} and references therein. This presents a severe challenge
for developing sufficiently valid analytic techniques for these equations, so numerical algorithms have been favored and
gradually emerged as essential alternatives in actual investigation.

The space-fractional PDEs constitute an important branch of PDEs in the studies on the evolution of complex dynamic behavior
governed by L\'{e}vy flights. As compared to the time-fractional PDEs, seeking the approximate solutions to such equations
appears to be more difficult by virtue of the vector structure of fractional Laplacian. Up to date, various numerical algorithms have
been designed, covering  finite difference methods \cite{Ref15,Ref14,R05,R07,R08,Ref13}, general Pad\'{e} approximation \cite{R04},
meshless point interpolation methods \cite{R09,R10}, moving least-squares meshless method \cite{xz01},
finite element methods \cite{Ref16,R12}, discontinuous Galerkin method (DGM) \cite{R14}, finite volume methods \cite{R15,R16},
spline approximation method (SAM) \cite{R17}, and spectral collocation methods \cite{Ref17,R18,R20}.
In \cite{R24,R21,R23,R22}, a series of operational matrix methods were constructed
based on the approximate expansions by using shifted Jacobi, Chebyshev, Legendre polynomials,
and Haar wavelets functions, as elements, respectively.
In \cite{Ref18}, a Crank--Nicolson alternative direction implicit Legendre-Galerkin spectral
method has been derived for nonlinear 2D Riesz space-fractional reaction-diffusion equations.
Du and Wang developed a finite element method (FEM) conjectured with a fast algebraic solver for the steady space-fractional diffusion equations 
by dividing the unit square domain into quadrilateral meshes \cite{Ref07}.
It is noteworthy that, nevertheless, all of the mentioned algorithms are only available for one-dimensional or
rectangular domain problems. Few works go for a valid method for the fractional problems on irregular domains.
Liu et. al handled a 2D fractional FitzHugh-Nagumo monodomain model by an implicit semi-alternative direction difference scheme
on circular domain and a union of half-circular and half-square domains with square meshes \cite{Ref06}, where the outer peripheral nodal points are in an approximate
match for boundaries. In \cite{Ref01}, we proposed a fully discrete FEM
for space-fractional Fisher's equation on rectangular domains and the generated meshes are relaxed to be unstructured; later, it
has been extended to the space-fractional diffusion equations on polygonal and elliptic domains \cite{Ref04}. Qiu et al.
developed a nodal DGM for the same type equations on unstructured triangular meshes and a L-shaped domain problem was
considered \cite{Ref05}. Pang et al. applied the RBFs meshless collocation methods
to solve the 2D space-fractional advection-diffusion equations on polygonal and  circular domains \cite{xz04}.
In general, finite difference method (FDM) is implemented on pre-defined meshes and inherits the shortcomings including the
difficulty in simulation of complex domain problems. DGM, FEM alleviate this issue, but both of them hinge on variational
principles and suffer huge computational burden and inflexibility to calculate the entities of their stiffness matrices that are on longer sparse.

Recent years have witnessed a keen interest in meshless techniques, 
 which avoid troublesome mesh generation and reconstruction,  or only use easily generable meshes in a flexible manner.
The meshless methods are superior to convectional mesh-dependent methods in terms of simulating complex flow, structure destruction, and extremely large deformation problems,
because they eliminate the element connectivity data and build trial functions entirely on scattered nodal points to decretize these problems.
Being a promising direction in computational mechanics, various meshless methods have been reported 
for numerical solutions of PDEs, typically including smooth particle hydrodynamics \cite{xz05}, finite point method \cite{xz06},
reproducing kernel particle method \cite{xz08}, partition of unity method \cite{xz13},
diffuse element method \cite{xz10}, finite sphere method,
element-free Galerkin methods \cite{xz02,xz03},
meshless local Petrov-Galerkin  method \cite{xz09}, boundary node method,
point interpolation method \cite{xz14}, and the others. We refer the readers to \cite{xz12,xz11} for their overall views. 

DQ method is understood as a direct numerical approach for finding the solutions of PDEs by reducing the equations
into ODEs via approximating the space derivatives as the weighted linear sums of the functional values at finite nodal points on  problem
domains. It was pioneered by Bellman and Casti \cite{R30}, and further discussed based on different basis functions, such as Lagrangian interpolation basis
functions, RBFs, orthogonal polynomials, Sinc and spline basis functions \cite{Ref20,Ref19,R37,R39,R36}. DQ method can achieve high accuracy by using a few nodal points;
besides, it is straight forward to implement and truly meshless.
In this context, regarding the current interest in fractional PDEs as effective models, we attempt to establish
valid DQ methods to solve the 2D space-fractional diffusion equations on arbitrary domains. The Multiquadric, Inverse Multiquadric, and
Gaussian RBFs are utilized as trial functions to determine the weighted coefficients
that we require to evaluate the fractional directional derivatives as the weighted linear sums of the functional values at regularly
distributed or scattered nodal points on problem domains; the Crank-Nicolson scheme is employed to advance the solutions in time. The proposed methods
extend the traditional DQ methods while inheriting their principal features.
The convergent behaviors of these techniques are studied on several numerical benchmarks with
a varying nodal number, including the square, trapezoidal, circular, and L-shaped domain problems.

The skeleton is organized as follows. In Section \ref{s2}, the space-fractional diffusion equations on a 2D domain are introduced for preliminaries.
Section \ref{s3} studies the DQ formulations for fractional directional derivatives and the way to determine their weights
by means of RBFs. In Section \ref{s4}, we construct time-stepping DQ methods by using the Crank-Nicolson scheme in time
and show the details of implementation. In Section \ref{s5}, illustrative tests on regular and irregular domains
are carried out to examine their accuracy and effectiveness. A concise remark is finally drawn in the last section. 

\section{Model problems}\label{s2}
\begin{figure*}
\centering
\includegraphics[width=3.8in]{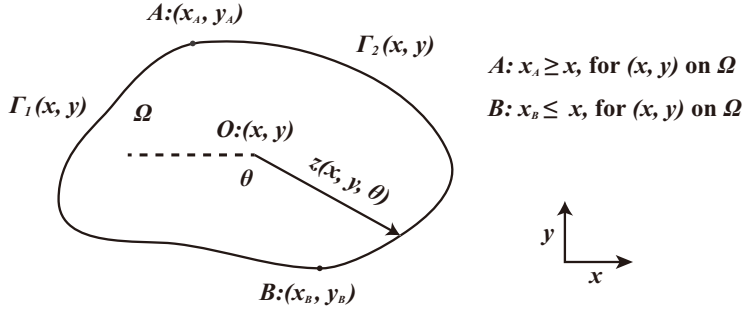}
\caption{The subsets of boundary $\Gamma_1(x,y)$, $\Gamma_2(x,y)$ and the distance $z(x,y,\theta)$ of $(x,y)$.}\label{fig0}
\end{figure*}

On a bounded domain $\Omega\subset \mathbb{R}^2$ with its boundary being $\partial\Omega$, the continuous
2D space-fractional diffusion equations are as follows
\begin{align}
\frac{\partial u(x,y,t)}{\partial t}-\kappa\!\int_{0}^{2\pi}\mathcal{D}^{\alpha}_{\theta}u(x,y,t)P(\theta)d\theta\!=\!f(x,y,t), \ (x,y;t)\in\Omega\times(0,T], \label{eq01}
\end{align}
subjected to the initial and boundary conditions
\begin{gather}
   u(x,y,0)=u_0(x,y),\quad (x,y)\in\Omega, \label{eq02}\\
    u(x,y,t)=g(x,y,t), \quad (x,y;t)\in\mathbb{R}^2\backslash \Omega\times(0,T], \label{eq03}
\end{gather}
where $1<\alpha\leqslant2$, $\kappa$ is the non-negative diffusivity parameter, $P(\theta)d\theta$
is a probability measure on the unit disk of $\mathbb{R}^2$,
and $f(x,y,t)$ is the pre-prescribed source function. In Eq. (\ref{eq01}), $\mathcal{D}^{\alpha}_{\theta}u(x,y,t)$
stands for the space-fractional directional derivatives defined in Caputo sense, i.e.,
\begin{align}\label{eq32}
    \mathcal{D}^{\alpha}_{\theta}u(x,y,t)=\mathcal{I}^{2-\alpha}_\theta \mathcal{D}^2_{\theta}u(x,y,t),\quad 0\leqslant\theta<2\pi,
\end{align}
with the integer-order directional derivatives
\begin{align*}
    \mathcal{D}^m_{\theta}v(x,y)=\bigg(\cos\theta\frac{\partial}{\partial x}+\sin\theta\frac{\partial}{\partial y}\bigg)^mv(x,y), \quad m\in \mathbb{N},
\end{align*}
and the fractional directional integration
\begin{align*}
    \mathcal{I}^{\mu}_\theta v(x,y)= \frac{1}{\Gamma(\mu)}\int_0^{z(x,y,\theta)}\omega^{\mu-1}v(x-\omega\cos\theta,y-\omega\sin\theta)d\omega,
\end{align*}
for $\mu>0$ while $\mathcal{I}^{\mu}_\theta v(x,y)=v(x,y)$ for $\mu=0$. $z(x,y,\theta)$ denotes the distance from the nodal point $(x,y)$ to $\partial\Omega$
in the direction $(-\cos\theta,-\sin\theta)$; see Fig. \ref{fig0}. 
In particular, when $\theta$ is taken to be some special values as $0$, $\pi$, $\mathcal{D}^{\alpha}_{\theta}u(x,y,t)$ recovers to
the commonly used Caputo derivatives with regard to $x$ \cite{Ref09}, i.e.,
\begin{align}
 \mathcal{D}^{\alpha}_{0}u(x,y,t)&=\frac{1}{\Gamma(2-\alpha)}
    \int^x_{\Gamma_1(x,y)}\frac{\partial^2 u(\omega,y,t)}{\partial \omega^2}\frac{d\omega}{(x-\omega)^{\alpha-1}},\label{eq30}\\
 \mathcal{D}^{\alpha}_{\pi}u(x,y,t)&=\frac{1}{\Gamma(2-\alpha)}
    \int^{\Gamma_2(x,y)}_x\frac{\partial^2 u(\omega,y,t)}{\partial \omega^2}\frac{d\omega}{(\omega-x)^{\alpha-1}},\label{eq31}
\end{align}
and Eqs. (\ref{eq01})-(\ref{eq03}) degenerate into the space-fractional diffusion equations in coordinate 
forms, which have been the topics of intense research, where $\Gamma_1(x,y)$, $\Gamma_2(x,y)$ are the subsets of $\partial\Omega$,
and satisfy $\partial\Omega=\Gamma_1(x,y)\cup \Gamma_2(x,y)$ and $\Gamma_1(x,y)\cap \Gamma_2(x,y)=\varnothing$,
for which, $A$, $B$ are dividing points; see Fig. \ref{fig0}.

In the sequel, without loss of generality, we consider the multi-term 2D space-fractional diffusion equations
with variable coefficients
\begin{align}
\frac{\partial u(x,y,t)}{\partial t}-\sum_{l=1}^{L}\kappa_l(x,y)\mathcal{D}^{\alpha_l}_{\theta_l}u(x,y,t)=f(x,y,t), \ \ (x,y;t)\in\Omega\times(0,T], \label{eq04}
\end{align}
subjected to the initial and boundary conditions
\begin{gather}
   u(x,y,0)=u_0(x,y),\quad (x,y)\in\Omega, \label{eq05}\\
   u(x,y,t)=g(x,y,t), \quad (x,y;t)\in\mathbb{R}^2\backslash \Omega\times(0,T],  \label{eq06}
\end{gather}
where $1<\alpha_l\leqslant2$, $0\leqslant\theta_l<2\pi$, $L\in\mathbb{Z}^+$, $\kappa_l(x,y)$ are the non-negative diffusivity parameters
but do not  fulfill $\sum_{l=1}^{L}\kappa_l(x,y)\equiv0$, and $\mathcal{D}^{\alpha_l}_{\theta_l}u(x,y,t)$ are the  
Caputo directional derivatives. Due to the characteristics of L\'{e}vy  processes, the boundary constrains of space-fractional diffusion
equations as Eqs. (\ref{eq01}), (\ref{eq04}) should be specified on the complementary set of $\Omega$ \cite{xz16}, so that
model problems can be well-posed. In this context, we assume that (\ref{eq06}) is always homogeneous outside $\Omega$, but may not be homogeneous at its boundary,
which says that the diffusive particles are killed after they escape from $\Omega$.

\section{DQ approximations for fractional directional derivatives}\label{s3}
In this section, we derive valid DQ formulations for approximating the fractional directional derivatives and study how to
compute the weighted coefficients by using Multiquadric, Inverse Multiquadric, and Gaussian RBFs as trial functions. To start with, let
$\Omega\subset \mathbb{R}^2$ be a bounded domain and $\textbf{x}_i$, $i=0,1,\ldots,M$, be a sequence of scattered nodal points distributed on $\Omega$ with
$\textbf{x}_i=(x_i,y_i)$. Define a set of proper basis functions $\{\phi_k(\textbf{x})\}_{k=0}^M$ that are adequately smooth to
guarantee the existence of fractional directional derivatives. DQ formulation can be described as combining evaluations
of a function at finite nodal points on problem domain so as to get a direct approximation to its partial derivative with regard to a variable.
In general, we always interpolate the exact solution of a fractional PDE like Eqs. (\ref{eq01})-(\ref{eq03}) in the form
\begin{equation}\label{eq07}
    u(\textbf{x},t)\cong\sum_{k=0}^M\delta_k(t)\phi_k(\textbf{x}),
\end{equation}
then enforce it to satisfy main equation as well as the related restrictions
at collocation points to determine the unknowns $\{\delta_k(t)\}_{k=0}^M$. However, if we have
\begin{align}
    \mathcal{D}^{\alpha}_{\theta}\phi_k(\textbf{x}_i)=\sum\limits_{j=0}^M {\omega_{ij}^{(\alpha)}\phi_k(\textbf{x}_j)},\ \ i,k=0,1,\ldots M,\label{eq08}
\end{align}
and substitute Eq. (\ref{eq08}) into Eq. (\ref{eq07}) after acting the differential-integral  operator
$\mathcal{D}^{\alpha}_{\theta}$ on both sides of Eq. (\ref{eq07}), it suffices to show that 
\begin{small}
\begin{align}
    \mathcal{D}^{\alpha}_{\theta}u(\textbf{x}_i,t)\cong\sum_{k=0}^M\delta_k(t)\mathcal{D}^{\alpha}_{\theta}\phi_k(\textbf{x}_i)
        =\sum_{k=0}^M\delta_k(t)\sum\limits_{j=0}^M {\omega_{ij}^{(\alpha)}\phi_k(\textbf{x}_j)}\cong \sum\limits_{j=0}^M {\omega_{ij}^{(\alpha)}u(\textbf{x}_j,t)},\label{eq09}
\end{align}
\end{small}
in the light of the linearity of $\mathcal{D}^{\alpha}_{\theta}$. In other words, Eqs. (\ref{eq09}) are valid DQ approximations to
$\mathcal{D}^{\alpha}_{\theta}u(\textbf{x}_i,t)$, $i=0,1,\ldots,M$, as long as Eqs. (\ref{eq08}) are fulfilled.
As a result, an approximate solution of a fractional PDE can be found by solving a set of ODEs without knowing $\{\delta_k(t)\}_{k=0}^M$.
We term $\omega_{ij}^{(\alpha)}$, $i,j=0,1,\ldots,M$, by the weighted or DQ coefficients for fractional directional derivatives, which
can be computed via Eqs. (\ref{eq08}) by reforming them in matrix-vector forms a priori;
$\{\phi_k(\textbf{x})\}_{k=0}^M$ are referred to as trial or test functions that can typically be chosen by those basis functions mentioned before. 
When $\alpha\in \mathbb{Z}^+$ and $\theta=0$, $\pi$, Eqs. (\ref{eq09}) reduce into the DQ formulations for classical derivatives.

\subsection{Radial basis functions}
The RBFs are known as a family of spline functions which are constructed from the distance between an arbitrary nodal point  
$\textbf{x}$ and its center $\textbf{x}_k$, i.e., $\varphi(||\textbf{x}-\textbf{x}_k||)$, abbreviated to $\varphi_k(\textbf{x})$, with the Euclidean norm $||\cdot||$.
These types of functions provide a set of excellent interpolating bases to interpolate the multivariable scattered data on high-dimensional domains and
have been served as an efficient tool in setting up truly meshless numerical algorithms for PDEs by virtue of their independency on geometric complexity
and the potential spectral accuracy of their interpolations. Among multiple types of RBFs,  
three familiar kinds of RBFs will be focused on hereinafter, i.e.,
\begin{itemize}
  \item Multiquadric RBFs:  $\varphi_k(\textbf{x})=\sqrt{||\textbf{x}-\textbf{x}_k||^2+\epsilon^2}$,
  \item Inverse Multiquadric RBFs: $\varphi_k(\textbf{x})=\frac{1}{\sqrt{||\textbf{x}-\textbf{x}_k||^2+\epsilon^2}}$,
  \item Gaussian RBFs:  $\varphi_k(\textbf{x})=\exp(-\epsilon^2||\textbf{x}-\textbf{x}_k||^2)$,
\end{itemize}
where $k=0,1,\ldots,M$, $M\in \mathbb{Z}^+$, and $\epsilon$ is a user number called by shape parameter. Given a function $y(\textbf{x},t)$
defined on $\Omega\subset \mathbb{R}^d$, $d=1,2,3$, its interpolating approximation based on these trial
bases can be written as 
\begin{equation}\label{eq10}
    y(\textbf{x},t)\cong\sum_{k=0}^{M}\lambda_k(t)\varphi_k(\textbf{x})+\sum_{s=1}^{M^*}\lambda_{M+s}(t)p_s(\textbf{x}),
\end{equation}
with the unknowns $\{\lambda_k(t)\}_{k=0}^M$, $\{\lambda_{M+s}(t)\}_{s=1}^{M^*}$ yet to be determined by collocating Eq. (\ref{eq10})
at nodal points $\{\textbf{x}_i\}_{i=0}^{M}$. In order to make the generated algebraic system to be well-posed, the orthogonal conditions
\begin{equation}\label{eq26}
   \sum_{k=0}^{M}\lambda_k(t)p_{s}(\textbf{x}_k)=0, \quad M^*=\frac{(r+d-1)!}{d!(r-1)!},
\end{equation}
for $s=1,2,\ldots,M^*$, should be imposed, in which, $\{p_s(\textbf{x})\}_{s=1}^{M^*}$
denote the basis functions of the polynomial space of degree at most $r-1$ on $\mathbb{R}^d$.
It is worthy to note that the shape parameter $\epsilon$ should be properly assigned in practical computation because it has a significant impact on
the approximate power of the radial basis interpolation, so is a RBFs-based method.


\subsection{Determination of weighted coefficients}
Before proceeding, we notice that the aforementioned polynomials in the right side of Eq. (\ref{eq10}) are not always necessary since
the augmented interpolating matrix is strictly positive definite for Inverse Multiquadrics, Gaussians, and is conditionally positive definite
for Multiquadrics under $M^*\geqslant1$ \cite{Ref28}. Hence, for ease of computing the weighted coefficients,
we take $M^*=0$ for Inverse Multiquadrics, Gaussians, and $M^*=1$ for Multiquadrics, in which case, the solvability of
the algebraic problem resulting from Eqs. (\ref{eq10})-(\ref{eq26}) is guaranteed.
For Inverse Multiquadrics and Gaussians, the approximate solution can be expressed by
$u(\textbf{x},t)\cong\sum_{k=0}^{M}\lambda_k(t)\varphi_k(\textbf{x})$. According to Eqs. (\ref{eq07})-(\ref{eq09}),
we then know the identity related to the fractional operator $\mathcal{D}^{\alpha}_{\theta}$, i.e.,
\begin{align*}
    \mathcal{D}^{\alpha}_{\theta}\varphi_k(\textbf{x}_i)=\sum\limits_{j=0}^M {\omega_{ij}^{(\alpha)}\varphi_k(\textbf{x}_j)},\ \ i,k=0,1,\ldots M,
\end{align*}
obtained by directly replacing $\phi_k(\textbf{x})$ by $\varphi_k(\textbf{x})$ in Eqs. (\ref{eq08}), which further lead to
the below solvable matrix-vector system:
\begin{align}\label{eq27}
\left( \begin{array}{cccc}
\varphi_0(\textbf{x}_0)&\varphi_0(\textbf{x}_1) &\cdots &\varphi_0(\textbf{x}_M)\\
\varphi_1(\textbf{x}_0)&\varphi_1(\textbf{x}_1)&\cdots &\varphi_1(\textbf{x}_M) \\
 \vdots &\vdots &\ddots&\vdots \\
\varphi_M(\textbf{x}_0)&\varphi_M(\textbf{x}_1)&\cdots &\varphi_M(\textbf{x}_M)
\end{array} \right)
\left( \begin{array}{c}
\omega^{(\alpha)}_{i0}\\
\omega^{(\alpha)}_{i1}\\
\vdots \\
\omega^{(\alpha)}_{iM}
\end{array} \right)
=\left( \begin{array}{c}
\mathcal{D}_\theta^\alpha\varphi_0(\textbf{x}_i)\\
\mathcal{D}_\theta^\alpha\varphi_1(\textbf{x}_i)\\
\vdots \\
\mathcal{D}_\theta^\alpha\varphi_M(\textbf{x}_i)
\end{array} \right),
\end{align}
with $i=0,1,\ldots M$ and the coefficient matrix being fully positive definite. As for Multiquadrics,
the polynomial term must be retained so as to ensure the well-posedness of the algebraic problems. Remembering that $M^*=1$,
there hold $u(\textbf{x},t)\cong\sum_{k=0}^{M}\lambda_k(t)\varphi_k(\textbf{x})+\lambda_{M+1}(t)$ and
the additional condition $\sum_{k=0}^{M}\lambda_k(t)=0$. Integrate these two equations into
a unified formula to get
\begin{equation*}
    u(\textbf{x},t)\cong\sum_{k=1}^{M}\lambda_k(t)[\varphi_k(\textbf{x})-\varphi_0(\textbf{x})]+\lambda_{M+1}(t).
\end{equation*}
According to Eqs. (\ref{eq07})-(\ref{eq09}), one then has
\begin{align}
    \mathcal{D}^{\alpha}_{\theta}[\varphi_k(\textbf{x}_i)-\varphi_0(\textbf{x}_i)]
        =\sum_{j=0}^M\omega_{ij}^{(\alpha)}[\varphi_k(\textbf{x}_j)-\varphi_0(\textbf{x}_j)],   \ \ i=0,1,\ldots M,\label{eq13}
\end{align}
for $k=1,2,\ldots,M$, while for $k=0$, it suffices to show  
\begin{equation}\label{eq15}
    \sum_{j=0}^M\omega_{ij}^{(\alpha)}=0,
\end{equation}
owing to $\mathcal{D}^{\alpha}_{\theta}C=0$ with a constant $C$. Rearranging Eqs. (\ref{eq13})-(\ref{eq15}) in matrix-vector forms,
a series of linear system of equations are finally obtained:
\begin{equation}\label{eq28}
\textbf{A}\boldsymbol{\omega}^{(\alpha)}_i=\textbf{D}_\theta^\alpha\boldsymbol{\varphi}(\textbf{x}_i)
     -\textbf{D}_\theta^\alpha\boldsymbol{\varphi}_0(\textbf{x}_i),\quad i=0,1,\ldots,M,
\end{equation}
where
\begin{align*}
\textbf{A}=\left( \begin{array}{cccc}
1&1 &\cdots &1\\
\varphi_1(\textbf{x}_0)-\varphi_0(\textbf{x}_0)&\varphi_1(\textbf{x}_1)-\varphi_0(\textbf{x}_1) &\cdots &\varphi_1(\textbf{x}_M)-\varphi_0(\textbf{x}_M)  \\
 \vdots &\vdots &\ddots&\vdots \\
\varphi_M(\textbf{x}_0)-\varphi_0(\textbf{x}_0) &\varphi_M(\textbf{x}_1)-\varphi_0(\textbf{x}_1)&\cdots &\varphi_M(\textbf{x}_M)-\varphi_0(\textbf{x}_M)
\end{array} \right),
\end{align*}
\begin{align*}
\boldsymbol{\omega}^{(\alpha)}_i=\left( \begin{array}{c}
\omega^{(\alpha)}_{i0}\\
\omega^{(\alpha)}_{i1}\\
\vdots \\
\omega^{(\alpha)}_{iM}
\end{array} \right),\
\textbf{D}_\theta^\alpha\boldsymbol{\varphi}(\textbf{x}_i)=\left( \begin{array}{c}
\mathcal{D}_\theta^\alpha\varphi_0(\textbf{x}_i)\\
\mathcal{D}_\theta^\alpha\varphi_1(\textbf{x}_i)\\
\vdots \\
\mathcal{D}_\theta^\alpha\varphi_M(\textbf{x}_i)
\end{array} \right),\
\textbf{D}_\theta^\alpha\boldsymbol{\varphi}_0(\textbf{x}_i)=\left( \begin{array}{c}
\mathcal{D}_\theta^\alpha\varphi_0(\textbf{x}_i)\\
\mathcal{D}_\theta^\alpha\varphi_0(\textbf{x}_i)\\
\vdots \\
\mathcal{D}_\theta^\alpha\varphi_0(\textbf{x}_i)
\end{array} \right).
\end{align*}

The unknown $\boldsymbol{\omega}^{(\alpha)}_i$, $i=0,1,\ldots,M$, are the weighted vectors what we seek for. However, here arises
another question of how to compute each component of the right-hand vectors
$\textbf{D}_\theta^\alpha\boldsymbol{\varphi}(\textbf{x}_i)$, $\textbf{D}_\theta^\alpha\boldsymbol{\varphi}_0(\textbf{x}_i)$,
which are the keys to formulating Eqs. (\ref{eq27}), (\ref{eq28}). Generally, the explicit expressions for the fractional
derivatives of a function can be derived but limited to a very small part of functions, that is to say,
the explicit expressions are unavailable in most cases, hence the explicit formulas should not be anticipated here, neither are the
conventional numerical quadrature rules since the weakly singular integral kernels in fractional derivatives.
In the sequel, we evaluate the fractional directional derivatives $\textbf{D}_\theta^\alpha\boldsymbol{\varphi}(\textbf{x}_i)$,
$\textbf{D}_\theta^\alpha\boldsymbol{\varphi}_0(\textbf{x}_i)$ by using Gauss-Jacobi quadrature rules
after making suitable integral transformations. At first, consider the components of $\textbf{D}_\theta^\alpha\boldsymbol{\varphi}(\textbf{x}_i)$ with
the trial functions being Inverse Multiquadrics. In view of the definition of second-order directional derivative, one has
\begin{align}\label{eq88}
   \mathcal{D}^2_{\theta}\varphi_k(\textbf{x})\!=\!\cos^2\theta\frac{\partial^2\varphi_k(x,y)}{\partial x^2}
    \!+\sin^2\theta\frac{\partial^2\varphi_k(x,y)}{\partial y^2}\!+2\sin\theta\cos\theta
    \frac{\partial^2\varphi_k(x,y)}{\partial x\partial y},
\end{align}
with $0\leqslant\theta<2\pi$, $k=0,1,\ldots,M$, which results in 
\begin{align}
   \mathcal{D}^2_{\theta}\varphi_k(\textbf{x})&=\frac{3\cos^2\theta(x-x_k)^2+3\sin^2\theta(y-y_k)^2+6\sin\theta\cos\theta(x-x_k)(y-y_k)}{\big((x-x_k)^2+(y-y_k)^2+\epsilon^2\big)^{5/2}} \nonumber\\
     &\quad-\frac{1}{\big((x-x_k)^2+(y-y_k)^2+\epsilon^2\big)^{3/2}}. \label{eq29}
\end{align}
On the other hand, by the identity (\ref{eq32}), there exists
\begin{equation}\label{eq17}
   \mathcal{D}_\theta^\alpha\varphi_k(\textbf{x})= \frac{1}{\Gamma(2-\alpha)}\int_0^{z(x,y,\theta)}\omega^{1-\alpha}
    \mathcal{D}^2_{\theta}\varphi_k(x-\omega\cos\theta,y-\omega\sin\theta)d\omega,
\end{equation}
with $0\leqslant\theta<2\pi$ and $k=0,1,\ldots,M$. Doing the integral transformation
$\omega=0.5z(x,y,\theta)(1+\varsigma)$ in Eq. (\ref{eq17}) as it is done in \cite{Ref22}, it is easy to check
\begin{equation}\label{eq16}
   \mathcal{D}_\theta^\alpha\varphi_k(\textbf{x})= \frac{1}{\Gamma(2-\alpha)}\bigg(\frac{z(x,y,\theta)}{2}\bigg)^{2-\alpha}\int_{-1}^1(1+\varsigma)^{1-\alpha}
    \chi(x,y,\theta,\varsigma)d\varsigma,
\end{equation}
where
\begin{equation*}
    \chi(x,y,\theta,\varsigma)=\mathcal{D}^2_{\theta}\varphi_k\bigg(x-\cos\theta\frac{z(x,y,\theta)(1+\varsigma)}{2},y-\sin\theta\frac{z(x,y,\theta)(1+\varsigma)}{2}\bigg),
\end{equation*}
which can be computed by Gauss-Jacobi quadrature rule. Let $\{c_s\}_{s=1}^Q$, $\{w_s\}_{s=1}^Q$ be quadrature
points and weights. Then, from Eq. (\ref{eq29}), it follows that
\begin{equation}\label{eq18}
   \mathcal{D}_\theta^\alpha\varphi_k(\textbf{x}_i)\cong\frac{1}{\Gamma(2-\alpha)}\bigg(\frac{z(x,y,\theta)}{2}\bigg)^{2-\alpha}\sum_{s=1}^Qw_s\eta(x_i,y_i,\theta,c_s),
\end{equation}
where
\begin{align*}
\eta(x,y,\theta,c_s)=\frac{3\cos^2\theta r_x^2+3\sin^2\theta r_y^2+6\sin\theta\cos\theta r_xr_y}{\left(r_x^2+r_y^2+\epsilon^2\right)^{5/2}}
    -\frac{1}{\left(r_x^2+r_y^2+\epsilon^2\right)^{3/2}},
\end{align*}
with the following quantities
\vspace{-1mm}
\begin{align*}
    &r_x=x-\cos\theta\frac{z(x,y,\theta)(1+c_s)}{2}-x_k,\\
    &r_y=y-\sin\theta\frac{z(x,y,\theta)(1+c_s)}{2}-y_k.
\end{align*}

In particular, when $\theta=0$, $\pi$, Eq. (\ref{eq32}) gives the ones in $x$-axis. Doing the integral transformations $\omega=x-0.5(x-\Gamma_1(x,y))(1+\varsigma)$ and
$\omega=x+0.5(\Gamma_2(x,y)-x)(1+\varsigma)$ in Eqs. (\ref{eq30})-(\ref{eq31}), respectively, $\mathcal{D}^{\alpha}_{0}\varphi_k(\textbf{x})$,
$\mathcal{D}^{\alpha}_{\pi}\varphi_k(\textbf{x})$  become
\begin{align*}
 \mathcal{D}^{\alpha}_{0}\varphi_k(\textbf{x})&=\frac{1}{\Gamma(2-\alpha)}\bigg(\frac{x-a}{2}\bigg)^{2-\alpha}
    \int^1_{-1}(1+\varsigma)^{1-\alpha}\varphi^{(2)}_k\bigg(x-\frac{(x-a)(1+\varsigma)}{2},y\bigg)d\varsigma,\\
 \mathcal{D}^{\alpha}_{\pi}\varphi_k(\textbf{x})&=\frac{1}{\Gamma(2-\alpha)}\bigg(\frac{b-x}{2}\bigg)^{2-\alpha}
    \int^1_{-1}(1+\varsigma)^{1-\alpha}\varphi^{(2)}_k\bigg(x+\frac{(b-x)(1+\varsigma)}{2},y\bigg)d\varsigma,
\end{align*}
with $a=\Gamma_1(x,y)$ and $b=\Gamma_2(x,y)$. Then, the computational formulas for $\mathcal{D}^\alpha_{0}\varphi_k(\textbf{x}_i)$,
$\mathcal{D}^\alpha_{\pi}\varphi_k(\textbf{x}_i)$, $i,\ k=0,1,\ldots,M$, are obtained as below
\begin{small}
\begin{align}
   &\mathcal{D}_0^\alpha\varphi_k(\textbf{x}_i)\cong\frac{1}{\Gamma(2-\alpha)}\bigg(\frac{x_i-a_i}{2}\bigg)^{2-\alpha}\sum_{s=1}^Qw_s\eta_0\bigg(x_i-\frac{(x_i-a_i)(1+c_s)}{2},y_i\bigg),\label{eq33}\\
   &\mathcal{D}_\pi^\alpha\varphi_k(\textbf{x}_i)\cong\frac{1}{\Gamma(2-\alpha)}\bigg(\frac{b_i-x_i}{2}\bigg)^{2-\alpha}\sum_{s=1}^Qw_s\eta_0\bigg(x_i+\frac{(b_i-x_i)(1+c_s)}{2},y_i\bigg),\label{eq35}
\end{align}
\end{small}
where $a_i=\Gamma_1(x_i,y_i)$, $b_i=\Gamma_2(x_i,y_i)$, and
\begin{align*}
\eta_0(x,y)=\frac{3(x-x_k)^2}{\big((x-x_k)^2+(y-y_k)^2+\epsilon^2\big)^{5/2}}-\frac{1}{\big((x-x_k)^2+(y-y_k)^2+\epsilon^2\big)^{3/2}}.
\end{align*}

When the trial functions are Gaussians, in the same manner, we can derive a computational formula analogous to (\ref{eq18})
for $\mathcal{D}^{\alpha}_{\theta}\varphi_k(\textbf{x}_i)$, but with
\begin{align*}
\eta(x,y,\theta,c_s)=2\epsilon^2e^{-\epsilon^2r_x^2-\epsilon^2r_y^2}\big(2\epsilon^2\cos^2\theta r_x^2 + 2\epsilon^2\sin^2\theta r_y^2 + 4\epsilon^2\sin\theta\cos\theta r_xr_y-1 \big),
\end{align*}
where $r_x$, $r_y$ are defined above. If $\theta$ is taken to be 0, $\pi$, the formulas for the fractional
derivatives in $x$-axis as (\ref{eq33})-(\ref{eq35}) can be obtained, but with
\begin{align*}
\eta_0(x,y)&=2\epsilon^2e^{-\epsilon^2(x-x_k)^2-\epsilon^2(y-y_k)^2}\big(2\epsilon^2(x-x_k)^2-1 \big).
\end{align*}
As for Multiquadrics, refer to \cite{xz04} for reference. When $\alpha=2$, $\mathcal{D}^{\alpha}_{\theta}\varphi_k(\textbf{x})$ recovers to
Eq. (\ref{eq88}), whose values are trivial. Once getting $\textbf{D}_\theta^\alpha\boldsymbol{\varphi}(\textbf{x}_i)$, $\textbf{D}_\theta^\alpha\boldsymbol{\varphi}_0(\textbf{x}_i)$,
$i=0,1,\ldots,M$, the weighted vectors $\boldsymbol{\omega}_{i}^{(\alpha)}$ are determined by solving Eqs. (\ref{eq27}), (\ref{eq28}) for each nodal point $\textbf{x}_i$ and the
fractional directional derivatives are then eliminated from a fractional PDE by using Eqs. (\ref{eq09}) as replacements, thus we can obtain the approximate solutions by treating ODEs instead.

\section{Time-stepping DQ methods and implemental processes}\label{s4}
In this section, we propose time-stepping DQ methods to approximate the solutions of Eqs. (\ref{eq04})-(\ref{eq06})
and show the implemental procedures. Define a lattice on $[0,T]$ with equally spaced points $t_n=n\tau$, $\tau=T/N$, $N\in\mathbb{Z}^+$, and
denote by $\omega_{ij}^{(\alpha_l)}$, $i$, $j=0,1,\ldots,M$, the weighted coefficients of
$\mathcal{D}_{\theta_l}^{\alpha_l} u(\textbf{x}_i,t)$.  On inserting the weighted sums
(\ref{eq09}) into Eq. (\ref{eq04}), we have the first-order ODEs:
\begin{align*}
\frac{\partial u(\textbf{x}_i,t)}{\partial t}-\sum_{l=1}^{L}\kappa_l(\textbf{x}_i)\sum\limits_{j=0}^M {\omega_{ij}^{(\alpha_l)}u(\textbf{x}_j,t)}
     =f(\textbf{x}_i,t), \ \ i=0,1,\cdots,M.
\end{align*}
Discretizing the ODEs by the Crank-Nicolson scheme reaches to
\begin{align}\label{eq19}
\begin{aligned}
&u(\textbf{x}_i,t_n)-\frac{\tau}{2}\sum_{l=1}^{L}\kappa_l(\textbf{x}_i)\sum\limits_{j=0}^M {\omega_{ij}^{(\alpha_l)}u(\textbf{x}_j,t_n)}=u(\textbf{x}_i,t_{n-1})\\
 &\quad\ +\frac{\tau}{2}\sum_{l=1}^{L}\kappa_l(\textbf{x}_i)\sum\limits_{j=0}^M {\omega_{ij}^{(\alpha_l)}u(\textbf{x}_j,t_{n-1})}+
      \tau f\bigg(\textbf{x}_i,t_n-\frac{\tau}{2}\bigg),
\end{aligned}
\end{align}
associated with the boundary constraints
\begin{equation}\label{eq20}
    u(\textbf{x}_i,t_n)=g(\textbf{x}_i,t_n),\ \ \textrm{for all}\ \textbf{x}_i\in\partial\Omega,
\end{equation}
where $i=0,1,\cdots,M$. The initial state are directly got from Eq. (\ref{eq05}).
In order to avoid wordy expressions below, we adopt $\kappa_i^l=\kappa_l(\textbf{x}_i)$,
$u^{n}_i=u(\textbf{x}_i,t_{n})$, $g^{n}_i=g(\textbf{x}_i,t_{n})$, and $f^{n-1/2}_i=f(\textbf{x}_i,t_n-\tau/2)$. Rewriting
Eqs. (\ref{eq19})-(\ref{eq20}) in matrix-vector form, a fully discrete DQ scheme reads
\begin{equation}\label{eq21}
  \bigg(\textbf{I}-\frac{\tau}{2}\sum_{l=1}^L\boldsymbol{\kappa}_l\textbf{W}_l\bigg)\textbf{U}^{n}=
  \bigg(\textbf{I}+\frac{\tau}{2}\sum_{l=1}^L\boldsymbol{\kappa}_l\textbf{W}_l\bigg)\textbf{U}^{n-1}+\tau\textbf{F}^{n-1/2},
\end{equation}
associated with the boundary constraints
\begin{equation}\label{eq22}
    u(\textbf{x}_i,t_n)=g(\textbf{x}_i,t_n),\ \ \textrm{for all}\ \textbf{x}_i\in\partial\Omega,
\end{equation}
where \textbf{I} is the $(M+1)\times(M+1)$ identity matrix,
$\textbf{U}^{n}=[u^n_0,u^n_1,\cdots,u^n_{M}]^T$,
$\boldsymbol{\kappa}_l=\textrm{diag}(\kappa_0^l,\kappa_1^l,\cdots,\kappa_{M}^l)$,
and $\textbf{W}_l$, $\textbf{F}^{n-1/2}$ are given by
\begin{align*}
\textbf{W}_l=\left( \begin{array}{cccc}
\omega^{(\alpha_l)}_{00}&\omega^{(\alpha_l)}_{01} &\cdots &\omega^{(\alpha_l)}_{0,M}\\
\omega^{(\alpha_l)}_{10}&\omega^{(\alpha_l)}_{11} &\cdots &\omega^{(\alpha_l)}_{1,M}  \\
 \vdots &\vdots &\ddots&\vdots \\
\omega^{(\alpha_l)}_{M,0}&\omega^{(\alpha_l)}_{M,1}&\cdots &\omega^{(\alpha_l)}_{M,M}
\end{array} \right),\quad
\textbf{F}^{n-1/2}=\left( \begin{array}{c}
f^{n-1/2}_0\\
f^{n-1/2}_1\\
 \vdots \\
f^{n-1/2}_{M}
\end{array} \right).
\end{align*}

In what follows, we show a detailed algorithm on how to implement Eqs. (\ref{eq21})-(\ref{eq22}).
Our codes are written on Matlab platform and for the ease of exposition, some commands will be
used as notations if no ambiguity is possible, for instance, $\textbf{A}(k,l)$ means
extracting the element from the matrix $\textbf{A}$, which locates at Row $k$ and Column $l$, if
the corresponding element exists. Let $nonb$, $boup$ be the index of the internal and boundary
nodal points on $\Omega$, respectively. Also, let $\tilde{\textbf{U}}^{n}$ be the unknowns related to the internal nodal points, i.e.,
$\textbf{U}^{n}(nonb)$. Then, Eqs. (\ref{eq21})-(\ref{eq22}) can be integrated into a unified form
\begin{align}\label{eq23}
      \bigg(\tilde{\textbf{I}}-\frac{\tau}{2}\sum_{l=1}^L\tilde{\boldsymbol{\kappa}}_l\textbf{K}_l\bigg)\tilde{\textbf{U}}^{n}=
  \bigg(\tilde{\textbf{I}}+\frac{\tau}{2}\sum_{l=1}^L\tilde{\boldsymbol{\kappa}}_l\textbf{K}_l\bigg)\tilde{\textbf{U}}^{n-1}+\tau\textbf{H}^{n-1/2},
\end{align}
where $\tilde{ \textbf{I}}$ is the identity matrix of rank $length(nonb)$, $\tilde{\boldsymbol{\kappa}}_l=\boldsymbol{\kappa}_l(nonb,nonb)$,
$\textbf{K}_l=\textbf{W}_l(nonb,nonb)$, and the right-hand vector
\begin{equation}\label{eq24}
    \textbf{H}^{n-1/2}=\tilde{\textbf{F}}^{n-1/2}+\frac{1}{2}\sum_{p=1}^{P}\sum_{l=1}^L
     \tilde{\boldsymbol{\kappa}}_l\big(\tilde{\textbf{g}}^n(p)+\tilde{\textbf{g}}^{n-1}(p)\big)\textbf{G}_l(:,p),
\end{equation}
with $P=length(boup)$, $\tilde{\textbf{g}}^n=\textbf{g}^n(boup)$, $\tilde{\textbf{F}}^{n-1/2}=\textbf{F}^{n-1/2}(nonb)$,
$\textbf{G}_l=\textbf{W}_l(nonb,boup)$, and $\textbf{g}^{n}=[g^n_0,g^n_1,\cdots,g^n_{M}]^T$.
A detailed implementation of the algorithm for Eqs. (\ref{eq23})-(\ref{eq24}) is summarized in the following flowchart
\begin{itemize}
  \item [1.]\ Input $\alpha_l$, $\theta_l$, $\epsilon$, $M$, $N$, and allocate $\{t_n\}_{n=0}^{N}$, $\{\textbf{x}_i\}_{i=0}^{M}$
  \item [2.]\ Define two zero matrices $\textbf{U}$, $\textbf{W}_l$ and form $\textbf{A}$, $\boldsymbol{\kappa}_l$ 
  \item [3.]\ while $i\leqslant M$ do
  \item [4.]\ \quad \   identify the distance $z(x_i,y_i,\theta_l)$ for each $\textbf{x}_i$
  \item [5.]\ \quad \   compute $\textbf{D}_{\theta_l}^{\alpha_l}\boldsymbol{\varphi}(\textbf{x}_i)$, $\textbf{D}_{\theta_l}^{\alpha_l}\boldsymbol{\varphi}_0(\textbf{x}_i)$
                       by the formulas as (\ref{eq18})-(\ref{eq35}) for each $\textbf{x}_i$
  \item [6.]\ \quad \   solve  Eqs. (\ref{eq27}), (\ref{eq28}) to obtain $\boldsymbol{\omega}^{(\alpha_l)}_i$ and set $\textbf{W}_l(i+1,:)=(\boldsymbol{\omega}^{(\alpha_l)}_i)^T$
  \item [7.]\ end while
  \item [8.]\ Let $\tilde{\boldsymbol{\kappa}}_l=\boldsymbol{\kappa}_l(nonb,nonb)$, $\textbf{K}_l=\textbf{W}_l(nonb,nonb)$, and $\textbf{G}_l=\textbf{W}_l(nonb,boup)$
  \item [9.]\ while $n\leqslant N$ do
  \item [10.]\ \quad \ compute $\tilde{\textbf{F}}^{n-1/2}$ and set $\textbf{H}^{n-1/2}=\tilde{\textbf{F}}^{n-1/2}$
  \item [11.]\ \quad \ compute $\textbf{g}^n$, $\textbf{g}^{n-1}$ and set $\tilde{\textbf{g}}^n=\textbf{g}^n(boup)$, $\tilde{\textbf{g}}^{n-1}=\textbf{g}^{n-1}(boup)$
  \item [12.]\ \quad \ while $p\leqslant length(boup)$ do
  \item [13.]\ \quad\quad \quad \ $\textbf{H}^{n-1/2}=\textbf{H}^{n-1/2}+\frac{1}{2}\sum_{l=1}^L
                                   \tilde{\boldsymbol{\kappa}}_l\big(\tilde{\textbf{g}}^n(p)+\tilde{\textbf{g}}^{n-1}(p)\big)\textbf{G}_l(:,p)$
  \item [14.]\ \quad \ end while
  \item [15.]\ \quad \ $\tilde{\textbf{U}}^{n}=\bigg(\bigg(\tilde{\textbf{I}}+\frac{\tau}{2}\sum_{l=1}^L\tilde{\boldsymbol{\kappa}}_l\textbf{K}_l\bigg)
                       \tilde{\textbf{U}}^{n-1}+\tau\textbf{H}^{n-1/2}\bigg)\bigg\backslash \bigg(\tilde{\textbf{I}}-\frac{\tau}{2}\sum_{l=1}^L\tilde{\boldsymbol{\kappa}}_l\textbf{K}_l\bigg)$
  \item [16.]\ \quad \ $\textbf{U}(nonb,n)=\tilde{\textbf{U}}^{n}$, $\textbf{U}(boup,n)=\tilde{\textbf{g}}^n$
  \item [17.]\ end while
  \item [18.]\ Output $\textbf{U}$ and terminate program
\end{itemize}

It is visible that our methods are truly free of troublesome mesh generation, so are the background cells in meshless Galerkin methods;
all the information need about the nodal points are their co-ordinates, which may make sense to treat the fractional equations as the models
under concern, especially for high-dimensional problems. Moreover, all involved RBFs are $C^\infty(\mathbb{R}^2)$ and this
is a sufficient condition that ensures the existence of $\mathcal{D}^\alpha_{\theta}\varphi_k(\textbf{x}_i)$, $i=0,1,\ldots,M$.

\section{Illustrative examples}\label{s5}
In this part, a couple of illustrative examples are carried out to gauge the practical performance of our algorithm, including two
1D problems and the 2D problems on the square, trapezoidal, circular, and L-shaped domains.
For simplicity, we abbreviate our methods to MQ-DQ, IM-DQ, and GA-DQ methods in the order of first appearance of the RBFs we use.
In the computation, the algorithm in \cite{Ref23} is applied to generate the Gauss-Jacobi
quadrature points and weights and 50 quadrature points and weights are preferred in calculating
the weighted coefficients. As to the shape parameter, how to choose it is still be open, although continued efforts have been devoted to theoretically
or numerically seek its optimal value in interpolation \cite{Ref25,Ref27,Ref24,Ref26,Ref29}, at which, the errors are minimum.
One thing we have to bear in mind is that its value would remarkably affect the accuracy of RBFs-based methods,
which is a principle working for DQ methods. Since the shape parameter should be adjusted with the nodal parameter $M$
so as to keep the condition numbers of interpolating matrices within a reasonable limit, 
we artificially select $\epsilon$ in 1D cases and apply $\epsilon=c^*/(M+1)^{0.25}$ in 2D cases as some works did with $c^*>0$.
During the entire computational processes, the errors are measured by
\begin{align*}
  e_2(\tau,M)=\sqrt{\frac{1}{M+1}\sum_{j=0}^M\Big|u_{j}^{n}-U_{j}^{n}\Big|^2},\quad
  e_\infty(\tau,M)=\max\limits_{0\leqslant j\leqslant M}\Big|u_{j}^n-U_{j}^n\Big|,
\end{align*}
and the corresponding convergent rates are computed by \cite{xz15}
\begin{align*}
  \textrm{Cov. rate}=\frac{d*\log_2\big(e_\nu(\tau,M_1)/e_\nu(\tau,M_2)\big)}{\log_2(M_2/M_1)},\quad \nu=2,\ \infty, 
\end{align*}
where $d$ is the dimension, $M_1$, $M_2$ are the nodal parameters as $M$ $(M_1\neq M_2)$, and
$u^n_{j}$, $U^n_{j}$ denote the exact and numerical solutions on the nodal point $\textbf{x}_j$ at the time level $n$, respectively.
For 1D problems, we perform the algorithm on the nodal distribution $x_j=0.5(1-\cos\frac{j\pi}{M})\ell+a$, $\ell=b-a$, $j=0,1,\cdots,M$,
on the computational interval $[a,b]$, while for 2D problems, both regular and irregular nodal distributions would be adopted in the tests. \\

\noindent
\textbf{Example 5.1.} Approximate the fractional derivative $\mathcal{D}_\pi^{\alpha}(1-x)^z$ on the interval $[0,1]$. From the
basic property of fractional derivative, we have
\begin{equation*}
   \mathcal{D}_\pi^{\alpha}(1-x)^z=\frac{\Gamma(z+1)}{\Gamma(z+1-\alpha)}(1-x)^{z-\alpha}.
\end{equation*}
In order to show the effectiveness of the DQ formulations, we place 11, 16, 21, and 26 nodal points on the computational interval, respectively,
and select the corresponding $\epsilon$'s for Multiquadrics by 0.3112, 0.2150, 0.1678, 0.1374, and for Inverse Multiquadrics by
0.4327, 0.3328, 0.2694, and 0.2255. As for Gaussians, we employ the data sets consisting of 4.0381, 5.3768, 6.6514, and 7.8994.
The numerical results for $\alpha=1.2$ and $z=3$ are tabulated in Table \ref{tab1}. Here, the approximations improve as $M$ increases, which implies that
the DQ formulations are valid. Besides, under these $\epsilon$'s, MQ-DQ formulation is more efficient than
IM-DQ and GA-DQ formulations in term of overall accuracy.\\

\begin{table*}
\centering
\caption{The numerical results with $\alpha=1.2$ and $z=3$ for Example 5.1}\label{tab1}
\begin{tabular}{lllllll}
\toprule
\multicolumn{1}{l}{\multirow{2}{0.5cm}{$M$}}
&\multicolumn{2}{l}{ MQ-DQ method} &\multicolumn{2}{l}{ IM-DQ method} &\multicolumn{2}{l}{GA-DQ method} \\
\cline{2-7}& $e_2(\tau,M)$  &$e_\infty(\tau,M)$ &$e_2(\tau,M)$  &$e_\infty(\tau,M)$ &$e_2(\tau,M)$  &$e_\infty(\tau,M)$ \\
\midrule 10    &2.5459e-02  &4.7254e-02  &  3.8207e-02  &6.2084e-02 &  9.3444e-02  &1.5869e-01 \\
         15    &9.8161e-03  &2.0683e-02  &  1.1916e-02  &2.5528e-02 &  3.2316e-02  &7.0755e-02 \\
        20     &4.8985e-03  &1.1519e-02  &  6.1154e-03  &1.3830e-02 &  1.6083e-02  &3.8576e-02\\
        25     &2.8489e-03  &7.1813e-03  &  3.5079e-03  &8.5431e-03 &  9.5893e-03  &2.4149e-02 \\
\bottomrule
\end{tabular}
\end{table*}

\begin{figure}[!htb]
\begin{minipage}[t]{0.49\linewidth}
\includegraphics[width=2.55in]{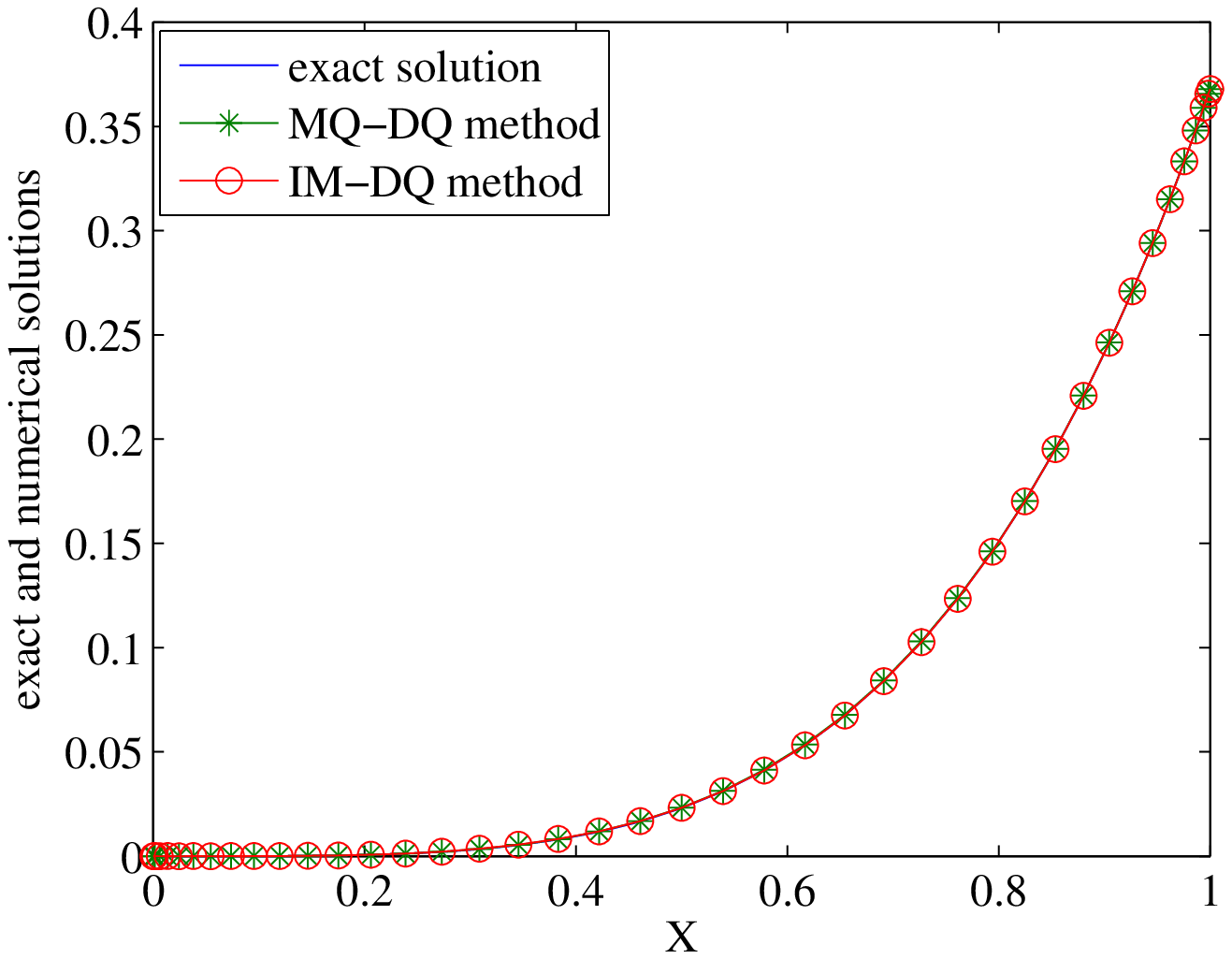}
\end{minipage}
\begin{minipage}[t]{0.49\linewidth}
\includegraphics[width=2.55in]{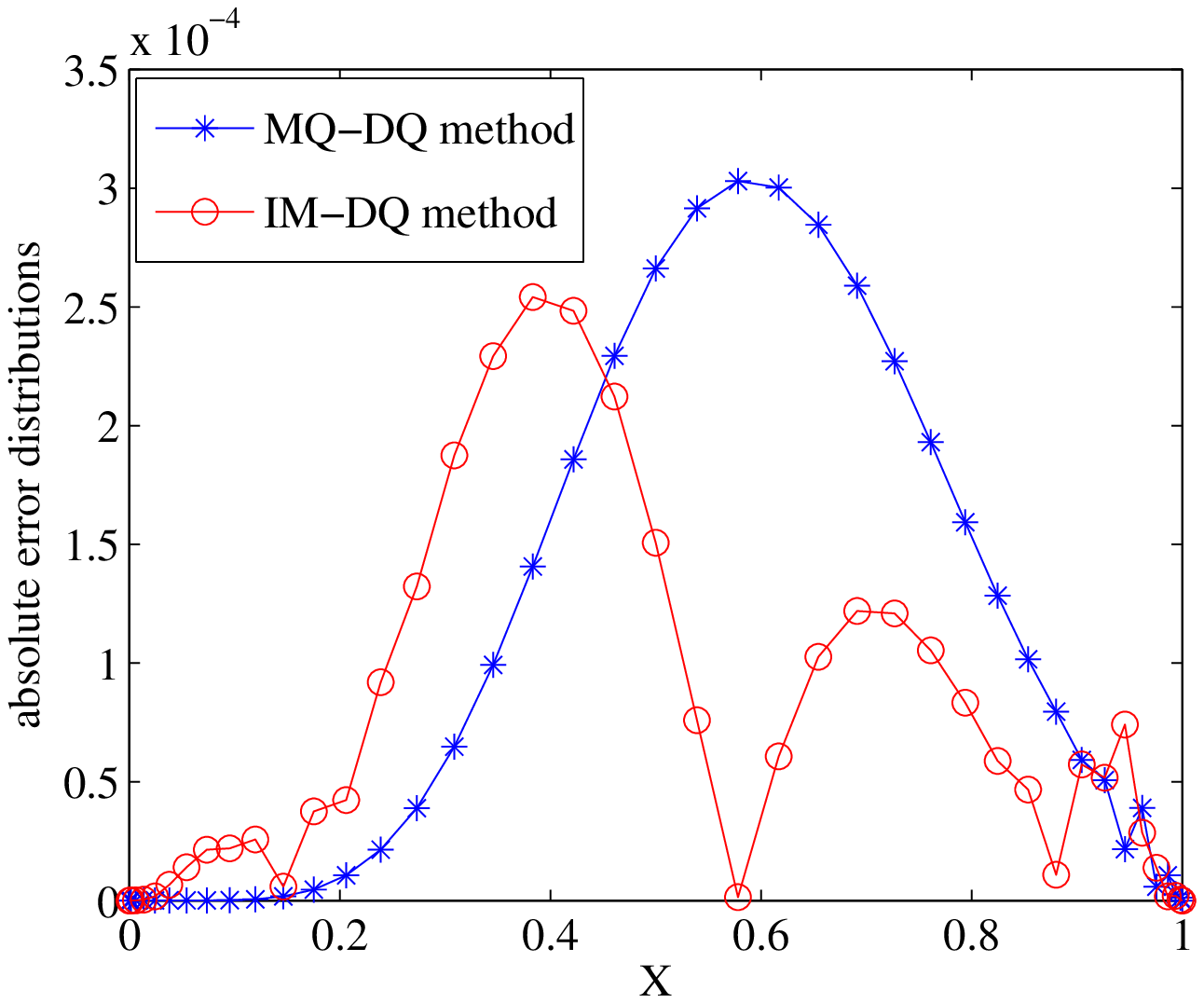}
\end{minipage}
\caption{The exact and numerical solutions (left) and absolute error distributions (right) at $t=1$
 by using $M=40$ and $\epsilon=0.0312$, $0.0511$ for MQ-DQ and IM-DQ methods, respectively.}\label{fig1}
\end{figure}

\noindent
\textbf{Example 5.2.} Letting $\theta=0$,  consider the fractional diffusion equation
\begin{align*}
\frac{\partial u(x,t)}{\partial t}-\frac{x^\alpha\Gamma(5-\alpha)}{24}\mathcal{D}^{\alpha}_{\theta}u(x,t)=-2e^{-t}x^4,
\end{align*}
on the interval $[0,1]$, subjected to the initial and boundary conditions
$u(x,0)=x^4$, $u(0,t)=0$, and $u(1,t)=e^{-t}$. It is verified that the exact solution is
$u(x,t)=e^{-t}x^4$. As in the first example, we place 16, 21, 26, and 31 nodal points on the
interval, respectively, and select $\epsilon$'s by 0.1875, 0.1128, 0.0712, 0.0613
for Multiquadrics and by 0.3098, 0.2135, 0.1567, and 0.1149 for Inverse Multiquadrics.
Taking $N=M$, the numerical results of MQ-DQ and IM-DQ methods at $t=1$ for $\alpha=1.5$ are compared with those obtained by SAM \cite{R17},
all of which are presented side by side in Table \ref{tab2}. Resetting $\epsilon=0.0312$ for Multiquadrics and $\epsilon=0.0511$ for Inverse Multiquadrics,
we display the exact and numerical solutions at $t=1$ with $M=40$ in Fig. \ref{fig1} (left)
and the absolute error distributions in Fig. \ref{fig1} (right). As seen from these table and graphs, DQ methods yield the approximations well
agreeing with the exact solutions and produce less errors than SAM; moreover, the accuracy of DQ methods clearly not only
depends on the nodal number, but also the shape parameter $\epsilon$. \\

\begin{table*}
\centering
\caption{A comparison of SAM and DQ methods at $t=1$ when $N=M$ and $\alpha=1.5$.}\label{tab2}
\begin{tabular}{llclclc}
\toprule
\multicolumn{1}{l}{\multirow{2}{0.5cm}{$M$}}
&\multicolumn{2}{l}{SAM \cite{R17}} &\multicolumn{2}{l}{ MQ-DQ method} &\multicolumn{2}{l}{IM-DQ method} \\
\cline{2-7}& $e_\infty(\tau,M)$  &Cov. rate &$e_\infty(\tau,M)$  &Cov. rate &$e_\infty(\tau,M)$  &Cov. rate \\
\midrule 15    &7.660e-04   &-    &  2.5379e-04  &-        &2.9346e-04  &- \\
         20    &4.493e-04   &1.9  &  1.3366e-04  &2.2288   &1.5818e-04  &2.1483\\
         25    &2.929e-04   &1.9  &  8.2231e-05  &2.1770   &9.8308e-05  &2.1315\\
         30    &2.067e-04   &1.9  &  5.5969e-05  &2.1102   &6.6635e-05  &2.1329\\
\bottomrule
\end{tabular}
\end{table*}

\begin{figure}[!htb]
\begin{minipage}[t]{0.49\linewidth}
\includegraphics[width=2.7in]{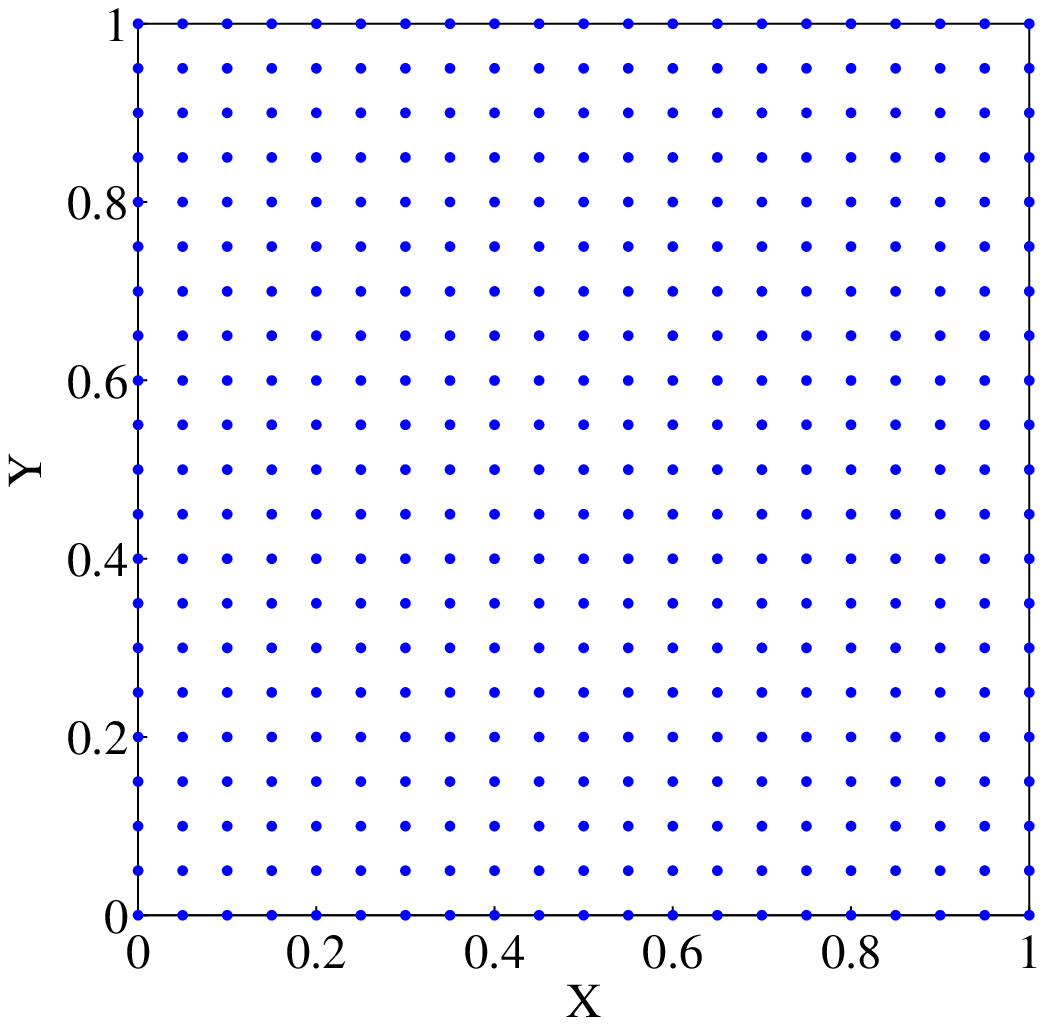}
\end{minipage}
\begin{minipage}[t]{0.49\linewidth}
\includegraphics[width=2.85in]{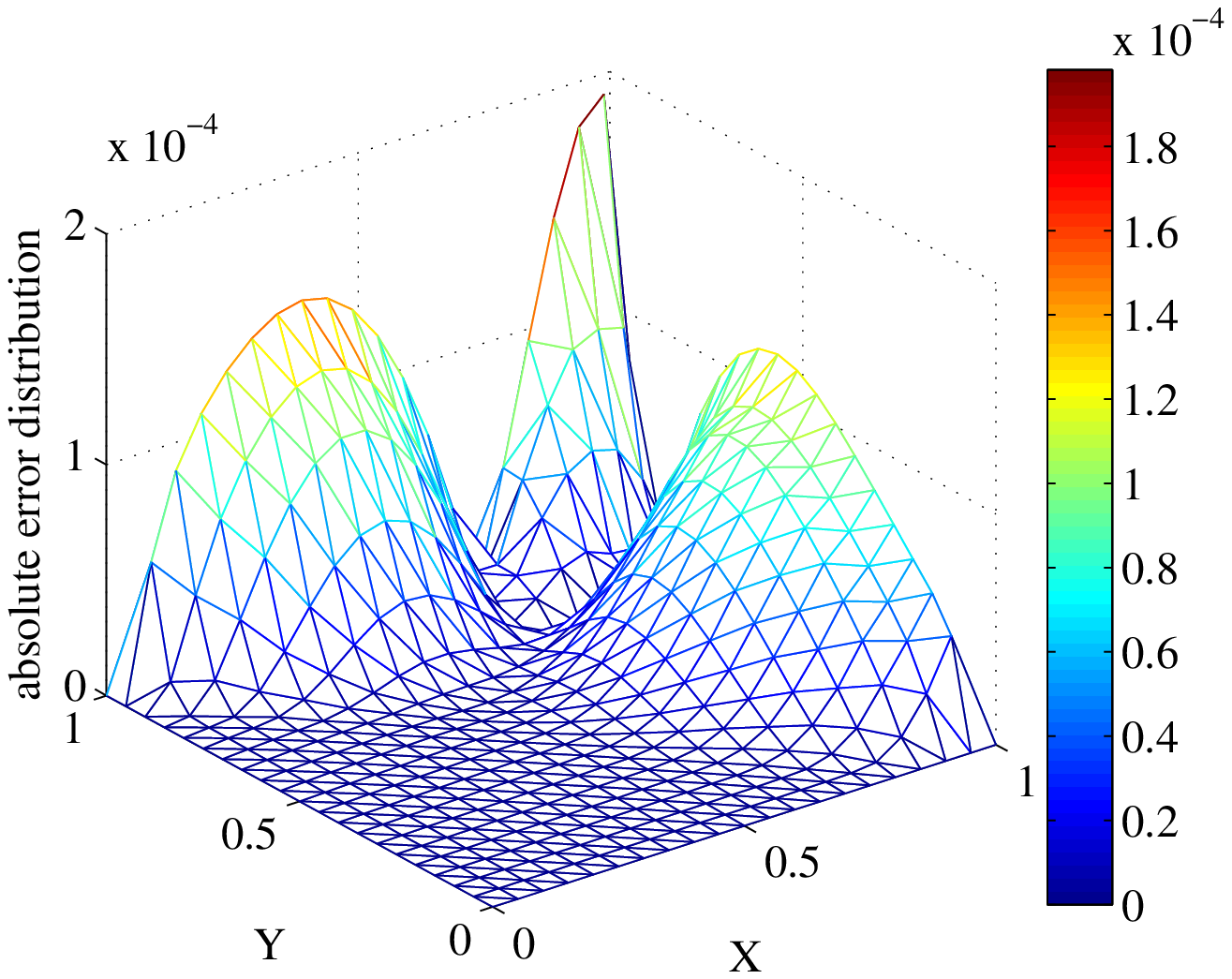}
\end{minipage}
\caption{The nodal distribution with $M=440$ and the absolute error distribution of MQ-DQ method.}\label{fig2}
\end{figure}

\begin{figure}[!htb]
\begin{minipage}[t]{0.48\linewidth}
\includegraphics[width=2.75in]{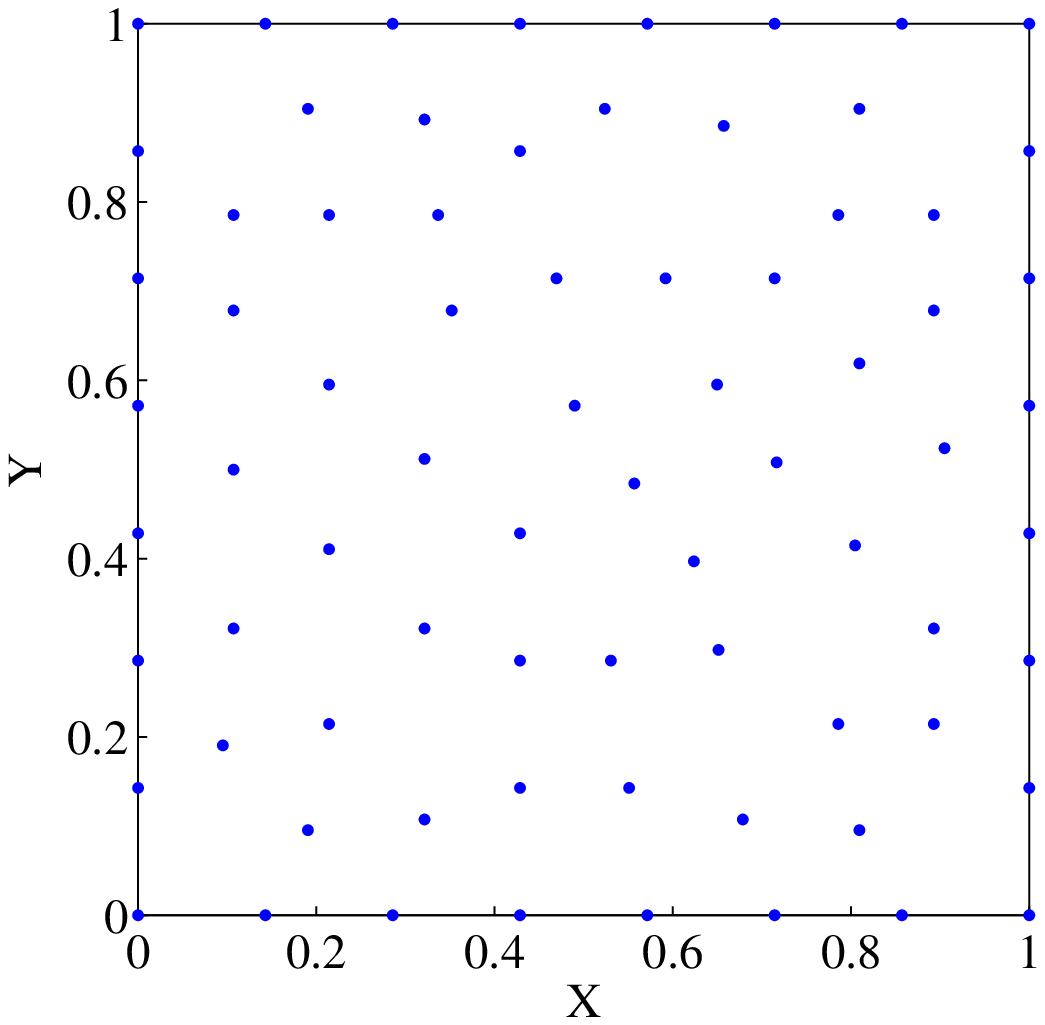}
\end{minipage}
\begin{minipage}[t]{0.48\linewidth}
\includegraphics[width=2.75in]{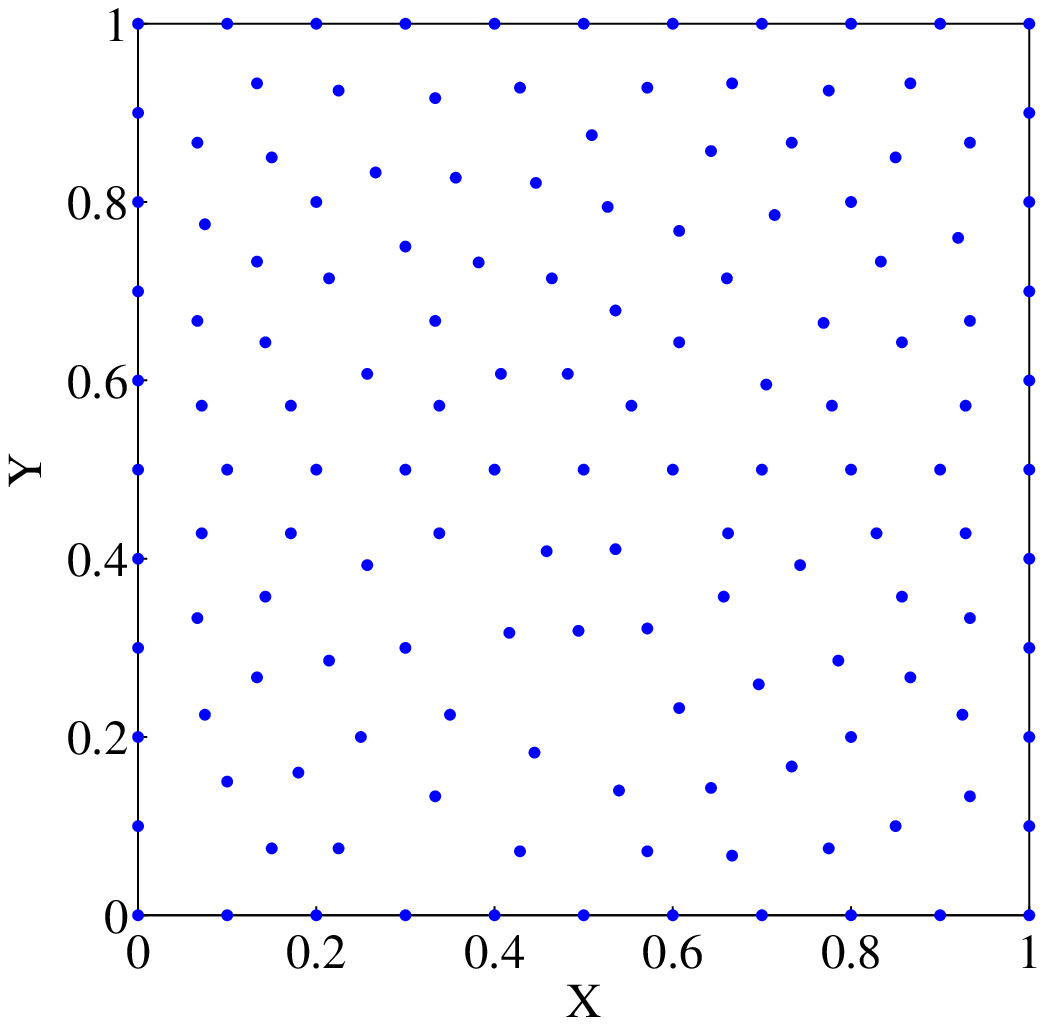}
\end{minipage}\\
\begin{minipage}[t]{0.48\linewidth}
\includegraphics[width=2.75in]{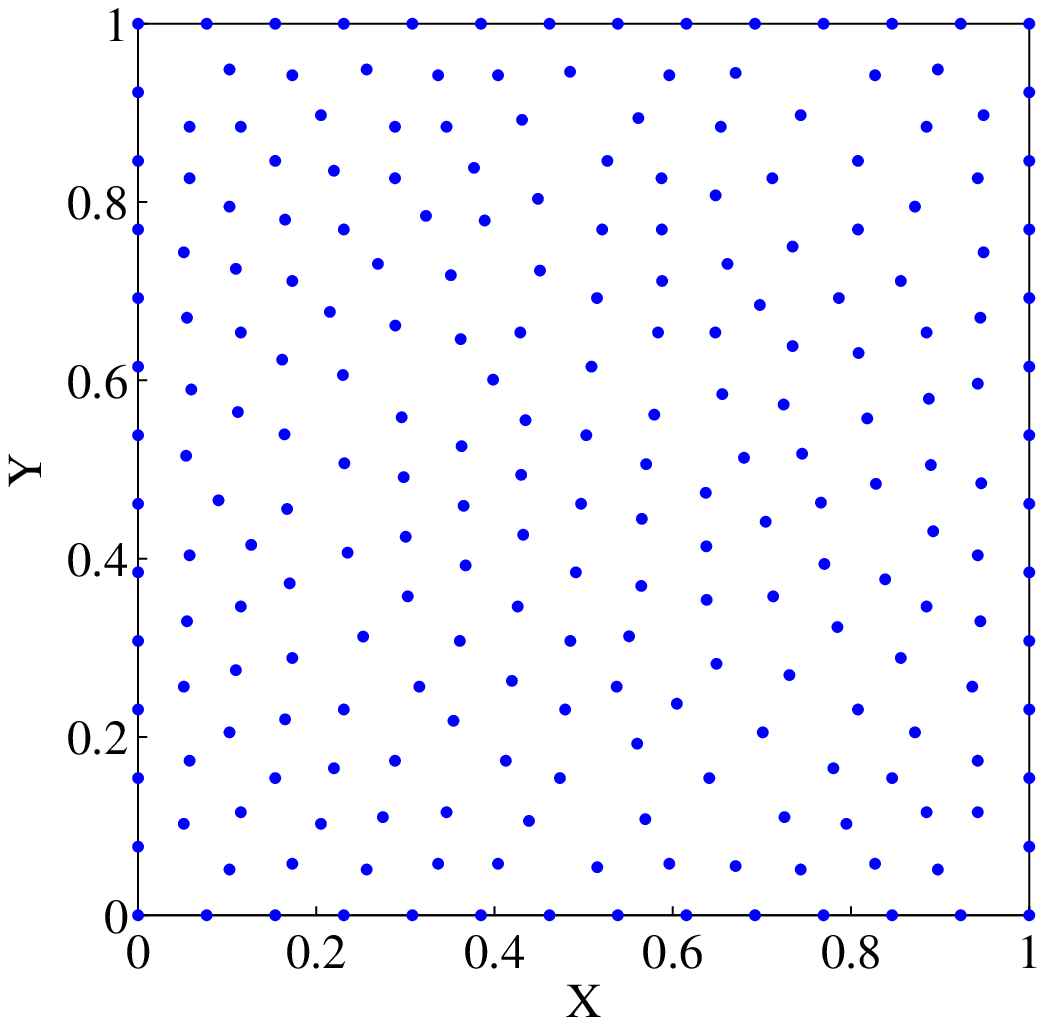}
\end{minipage}
\begin{minipage}[t]{0.48\linewidth}
\includegraphics[width=2.75in]{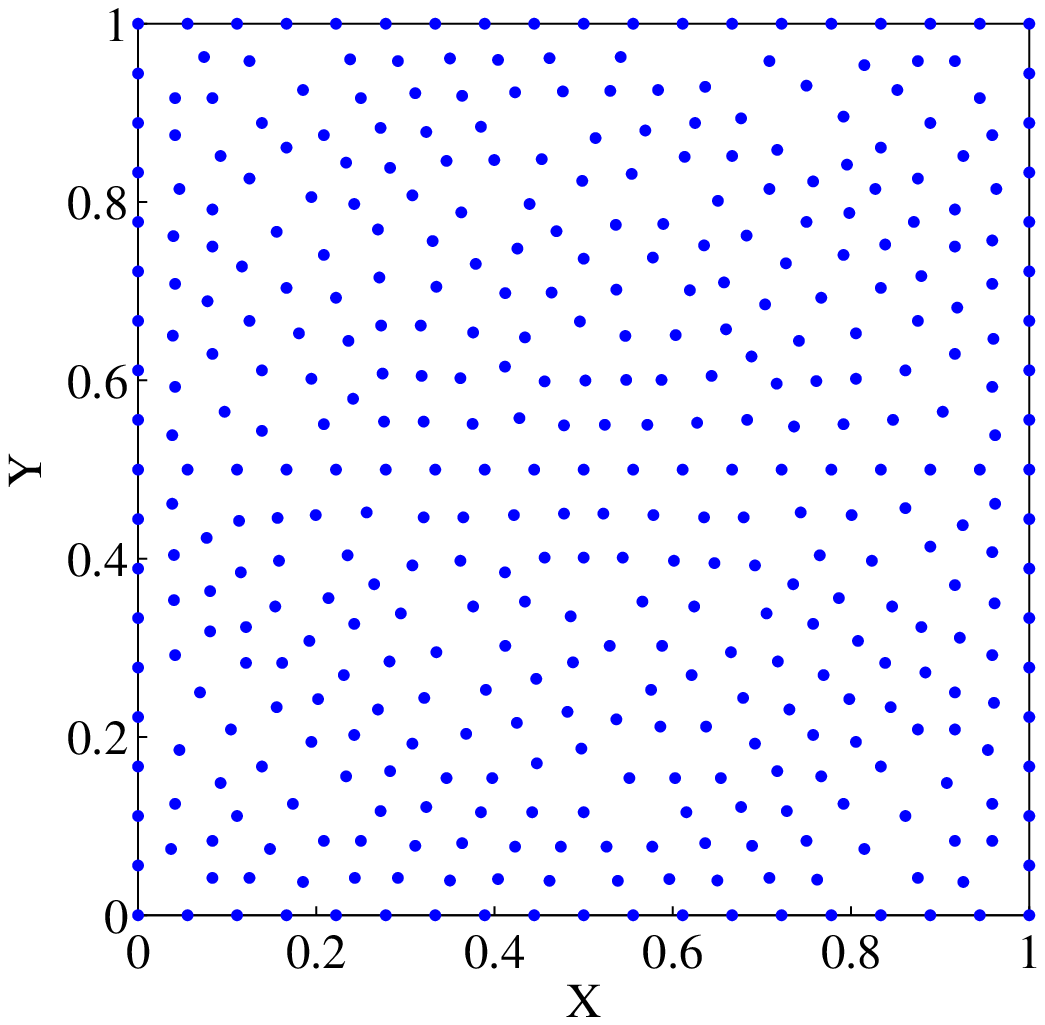}
\end{minipage}
\caption{The nodal distributions with nodal numbers $74$, $144$, $234$, and $424$ for Example 5.3.}\label{fig3}
\end{figure}

\begin{figure}[!htb]
\begin{minipage}[t]{0.49\linewidth}
\includegraphics[width=2.45in]{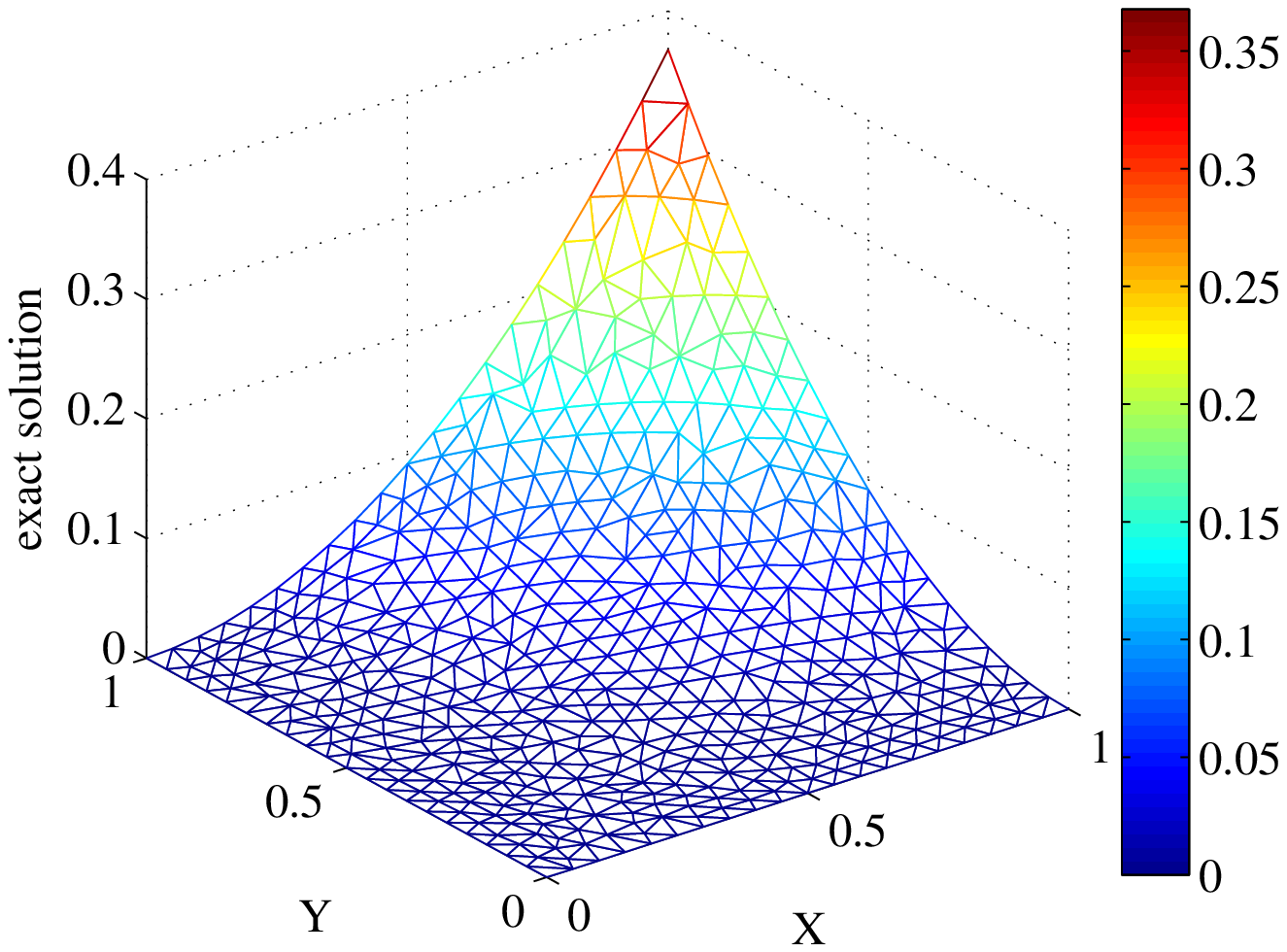}
\end{minipage}
\begin{minipage}[t]{0.49\linewidth}
\includegraphics[width=2.45in]{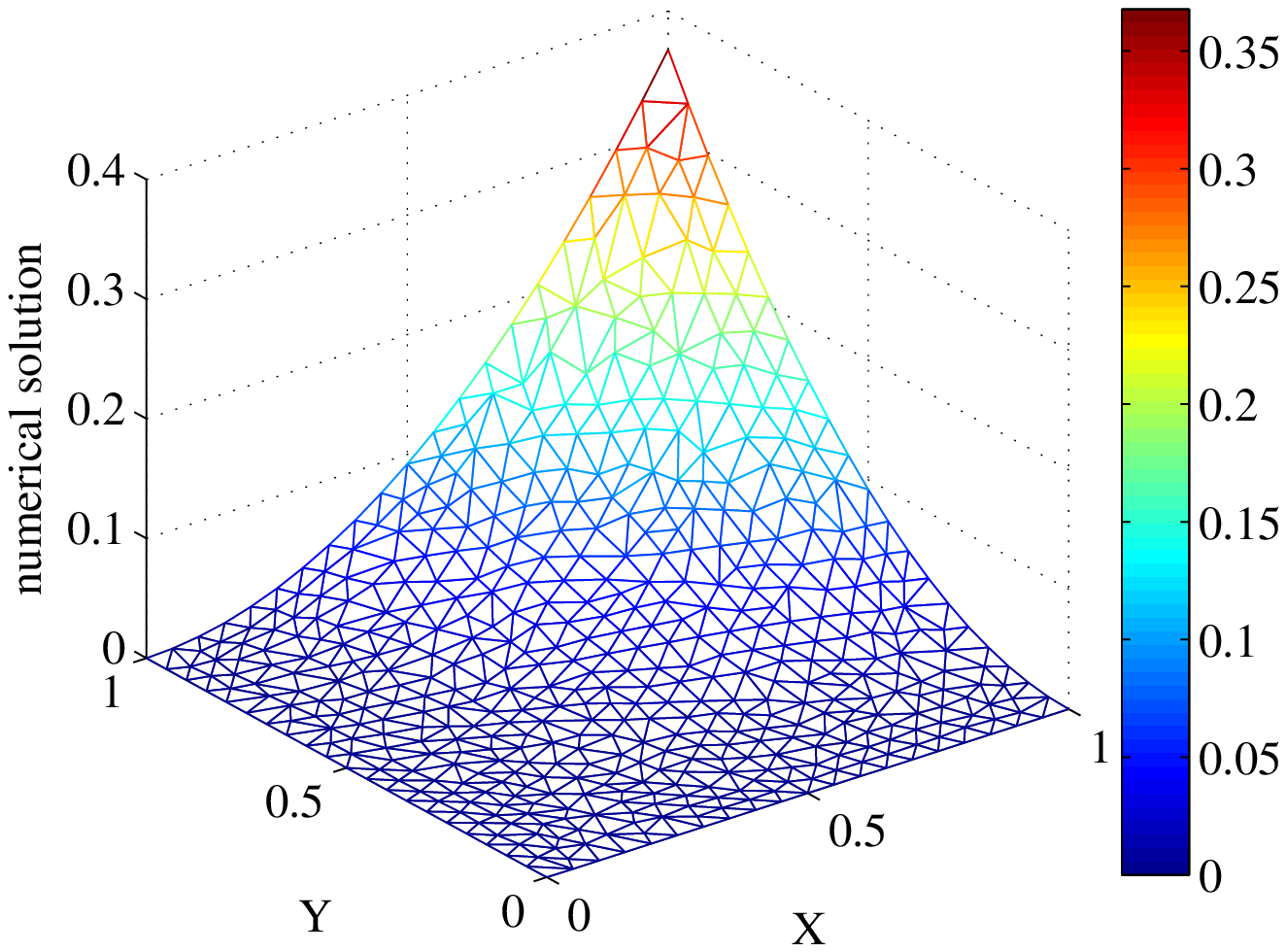}
\end{minipage}
\caption{The exact and numerical solutions created by GA-DQ method for Example 5.3.}\label{fig4}
\end{figure}

\noindent
\textbf{Example 5.3.} Consider the 2D fractional diffusion equation
\begin{align*}
\frac{\partial u(x,y,t)}{\partial t}-\kappa_1(x,y)\mathcal{D}^{\alpha_1}_{\theta_1}u(x,y,t)-\kappa_2(x,y)\mathcal{D}^{\alpha_2}_{\theta_2}u(x,y,t)=f(x,y,t), \label{eq25}
\end{align*}
on the square domain $[0,1]\times[0,1]$ in two separate cases. The former utilizes the regularly distributed nodal points,
while the latter chooses the irregularly distributed nodal points.  To compare with FDM \cite{Ref14}, we
extract the regular nodal distribution from the structured meshes generated by FreeFem++.
The computational parameters of these two cases are arranged as follows:
\begin{itemize}
  \item [(i)] Letting $\alpha_1=1.8$, $\alpha_2=1.6$, $\theta_1=0$, $\theta_2=\pi/2$, $\kappa_1(x,y)=\frac{\Gamma(2.2)x^{2.8}y}{6}$,
             $\kappa_2(x,y)=\frac{2xy^{2.6}}{\Gamma(4.6)}$, and the source function $f(x,y,t)=-e^{-t}(1+2xy)x^3y^{3.6}$,
             solve the above equation with the initial and boundary conditions $u_0(x,y)=x^3y^{3.6}$, $u(0,y,t)=u(x,0,t)=0$, $u(1,y,t)=e^{-t}y^{3.6}$,
             and $u(x,1,t)=e^{-t}x^{3}$. The exact solution is given by $u(x,y,t)=e^{-t}x^3y^{3.6}$.
  \item [(ii)] Letting $\alpha_1=1.8$, $\theta_1=\pi/4$, $\kappa_1(x,y)=x^{\alpha_1}$, $\kappa_2(x,y)=0$, and the source function
              $f(x,y,t)=-e^{-t}x^2y^2-e^{-t}x^{\alpha_1} f^*(x,y,t)$ with
                \begin{align*}
                 f^*(x,y,t)=\left\{
                \begin{array}{ll}
                   \frac{2^{1-\frac{\alpha_1 }{2}} y^{2-\alpha_1}\big((\alpha_1-4)(\alpha_1-3) x^2-2(\alpha_1-4)\alpha_1 xy+(\alpha_1-1)\alpha_1  y^2\big)}{\Gamma(5-\alpha_1)}, \ \ x\geqslant y,\\
                   \frac{2^{1-\frac{\alpha_1 }{2}} x^{2-\alpha_1}\big((\alpha_1-1)\alpha_1 x^2-2(\alpha_1-4)\alpha_1 xy+(\alpha_1-4)(\alpha_1 -3) y^2\big)}{\Gamma(5-\alpha_1)}, \ \ x<y,
                \end{array}
                \right.
                \end{align*}
             solve the above equation with the initial and boundary conditions $u_0(x,y)=x^2y^2$, $u(0,y,t)=u(x,0,t)=0$, $u(1,y,t)=e^{-t}y^{2}$,
             and $u(x,1,t)=e^{-t}x^{2}$. It can be verified that the exact solution is $u(x,y,t)=e^{-t}x^2y^2$.
\end{itemize}
In the second case, the source function suffers relatively complex mathematical structure since $\theta_1$ does not along the axis directions.
In Table \ref{tab3}, we give a comparison between FDM and DQ methods in term of $e_\infty(\tau,M)$ at $t=1$ for the case of (i) with $N=\sqrt{M+1}-1$,
$c^*=0.98$ for Multiquadrics and $c^*=1.22$ for Inverse Multiquadrics.
The regular distribution of nodal points with $M=440$ and the corresponding absolute error distribution of MQ-DQ method
are shown in Fig. \ref{fig2}. For the case of (ii), we list the numerical results at $t=1$ of DQ methods in Table \ref{tab4}
by using the scattered nodal points of total numbers $74$, $144$, $234$, and $424$, with $N=2000$, $c^*=0.89$ for Multiquadrics and $c^*=1.25$
for Inverse Multiquadrics. In Fig. \ref{fig3}, we exhibit the used irregular distributions of nodal points, and in Fig. \ref{fig4}, by resetting
$M=518$ and $\epsilon=8.6924$, we compare the exact and numerical solutions at $t=1$ created by GA-DQ method.
As expected, DQ methods work well on both regular and irregular nodal distributions and yield the approximations that are in perfect agreement with
the exact ones. It is also drawn that our methods achieve the errors in the same scale of magnitude as FDM with less nodal numbers. \\

\begin{table*}
\centering
\caption{A comparison of FDM and DQ methods at $t=1$ when $N=\sqrt{M+1}-1$.}\label{tab3}
\begin{tabular}{llcclclc}
\toprule
\multicolumn{1}{l}{\multirow{2}{0.4cm}{$M$}}
&\multicolumn{2}{l}{FDM \cite{Ref14}} & \multicolumn{1}{l}{\multirow{2}{0.3cm}{$M$}} &\multicolumn{2}{l}{ MQ-DQ method} &\multicolumn{2}{l}{IM-DQ method} \\
\cline{2-3}\cline{5-8}&$e_\infty(\tau,M)$  &Cov. rate &  &$e_\infty(\tau,M)$  &Cov. rate &$e_\infty(\tau,M)$  &Cov. rate \\
\midrule 120    &1.2629e-03   &-    & 99    &  1.2391e-03  &-        &3.2338e-03  &- \\
         440    &6.7325e-04   &1.88 & 195   &  5.3030e-04  &2.5222   &1.5975e-03  &2.0959\\
         1680   &3.4824e-04   &1.93 & 288   &  3.3018e-04  &2.4405   &9.9305e-04  &2.4487\\
         6560   &1.7660e-04   &1.97 & 440   &  1.9823e-04  &2.4144   &5.5787e-04  &2.7290\\
\bottomrule
\end{tabular} 
\end{table*}

\begin{table*}[!htb]
\centering
\caption{The numerical results at $t=1$ with $N=2000$ and $\alpha=1.8$ for Example 5.3.}\label{tab4}
\begin{tabular}{lclclc}
\toprule
\multicolumn{1}{l}{\multirow{2}{0.6cm}{$M$}} &\multicolumn{2}{l}{ MQ-DQ method} &\multicolumn{2}{l}{IM-DQ method} \\
\cline{2-5}& $e_2(\tau,M)$  &$e_\infty(\tau,M)$  &$e_2(\tau,M)$  & $e_\infty(\tau,M)$ \\
\midrule 73     &  4.5310e-04  &1.6351e-03   &9.5684e-04  &2.8237e-03 \\
         143    &  2.6755e-04  &9.3234e-04   &4.6647e-04  &1.4306e-03 \\
         233    &  8.8580e-05  &3.5897e-04   &1.6103e-04  &6.1178e-04 \\
         423    &  2.6730e-05  &1.2532e-04   &4.7417e-05  &1.7907e-04 \\
\bottomrule
\end{tabular}
\end{table*}

\begin{figure}
\begin{minipage}[t]{0.49\linewidth}
\includegraphics[width=2.51in]{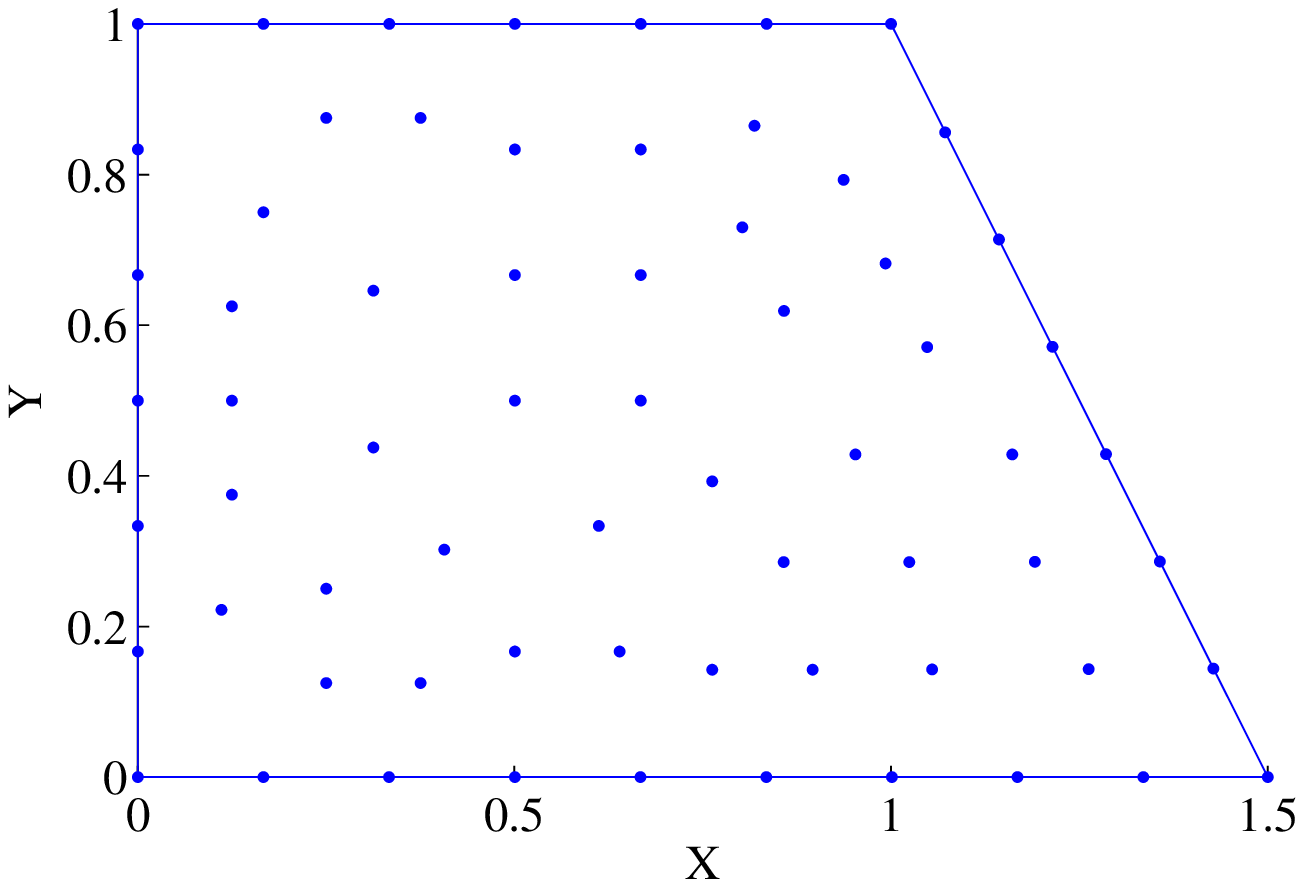}
\end{minipage}
\begin{minipage}[t]{0.49\linewidth}
\includegraphics[width=2.51in]{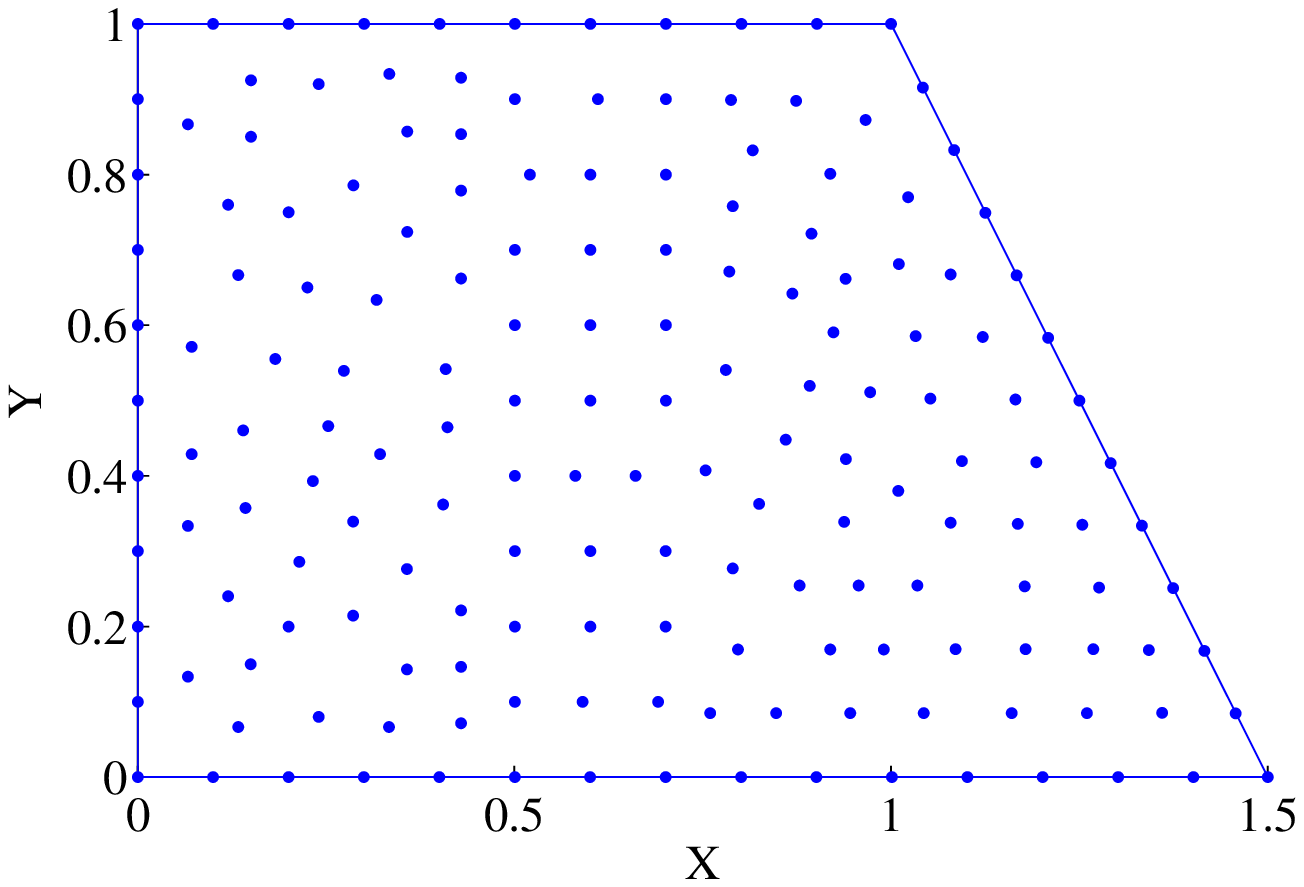}
\end{minipage}\\
\begin{minipage}[t]{0.49\linewidth}
\includegraphics[width=2.51in]{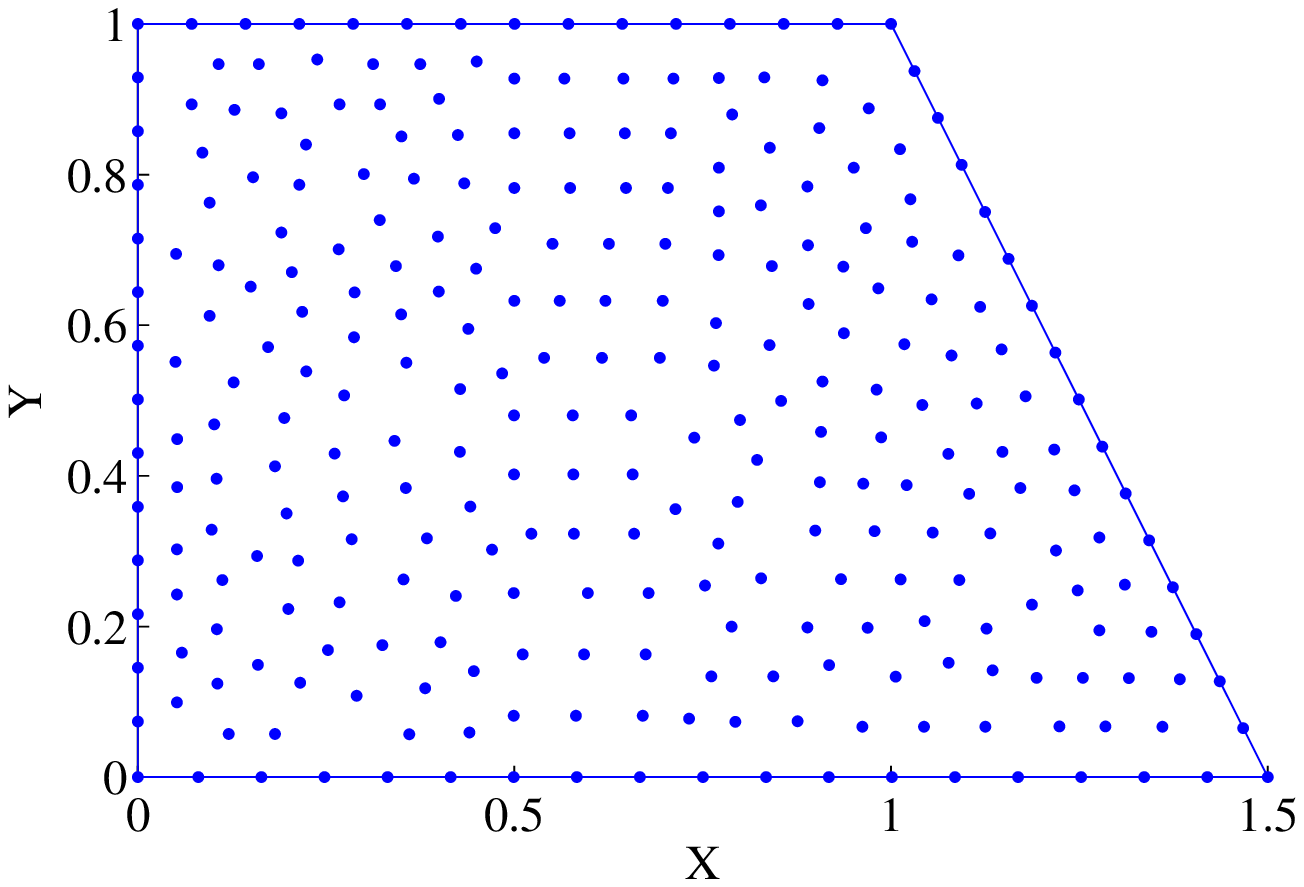}
\end{minipage}
\begin{minipage}[t]{0.49\linewidth}
\includegraphics[width=2.51in]{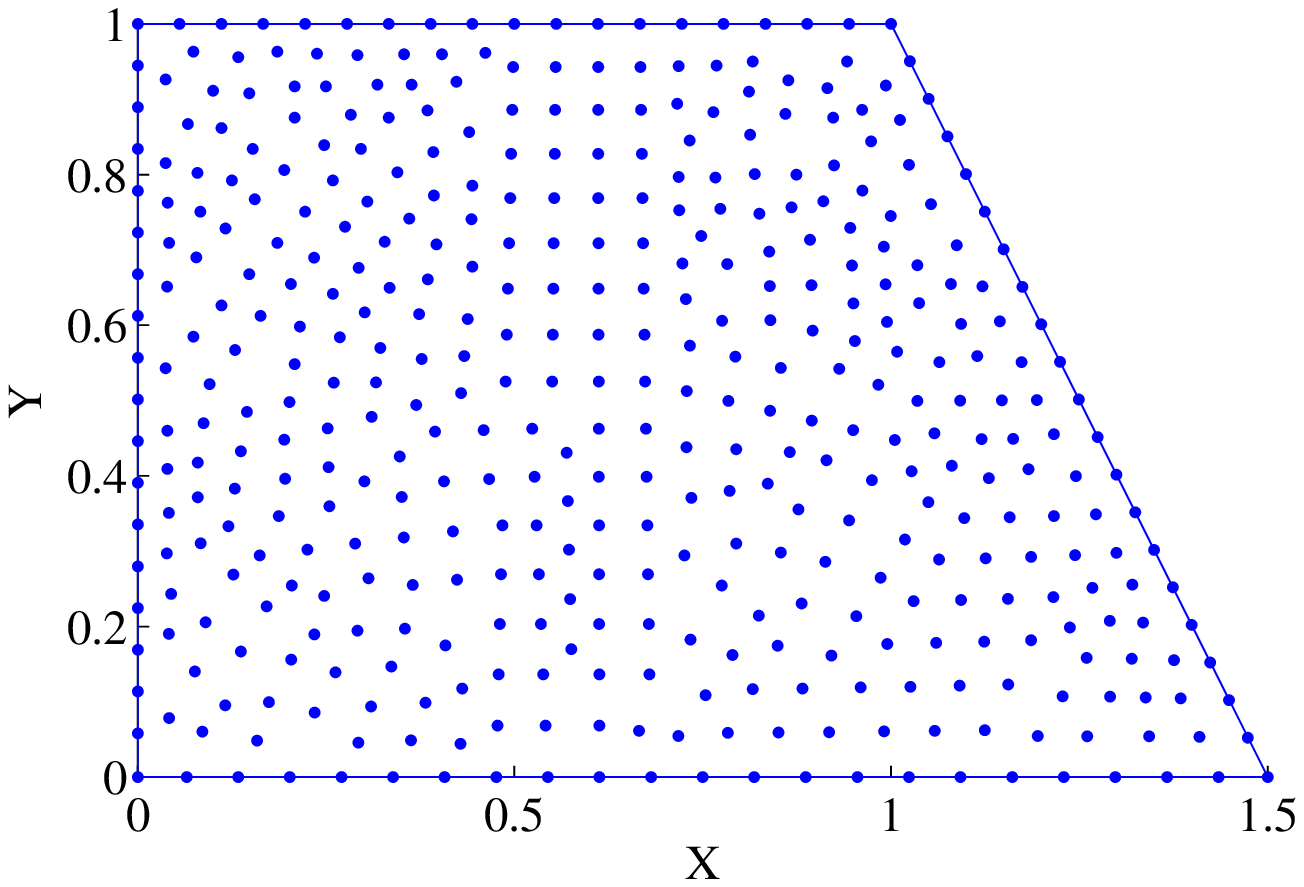}
\end{minipage}
\caption{The nodal distributions with nodal numbers $66$, $171$, $287$, and $437$ for Example 5.4.}\label{fig5}
\end{figure}

\begin{figure}
\begin{minipage}[t]{0.49\linewidth}
\includegraphics[width=2.45in]{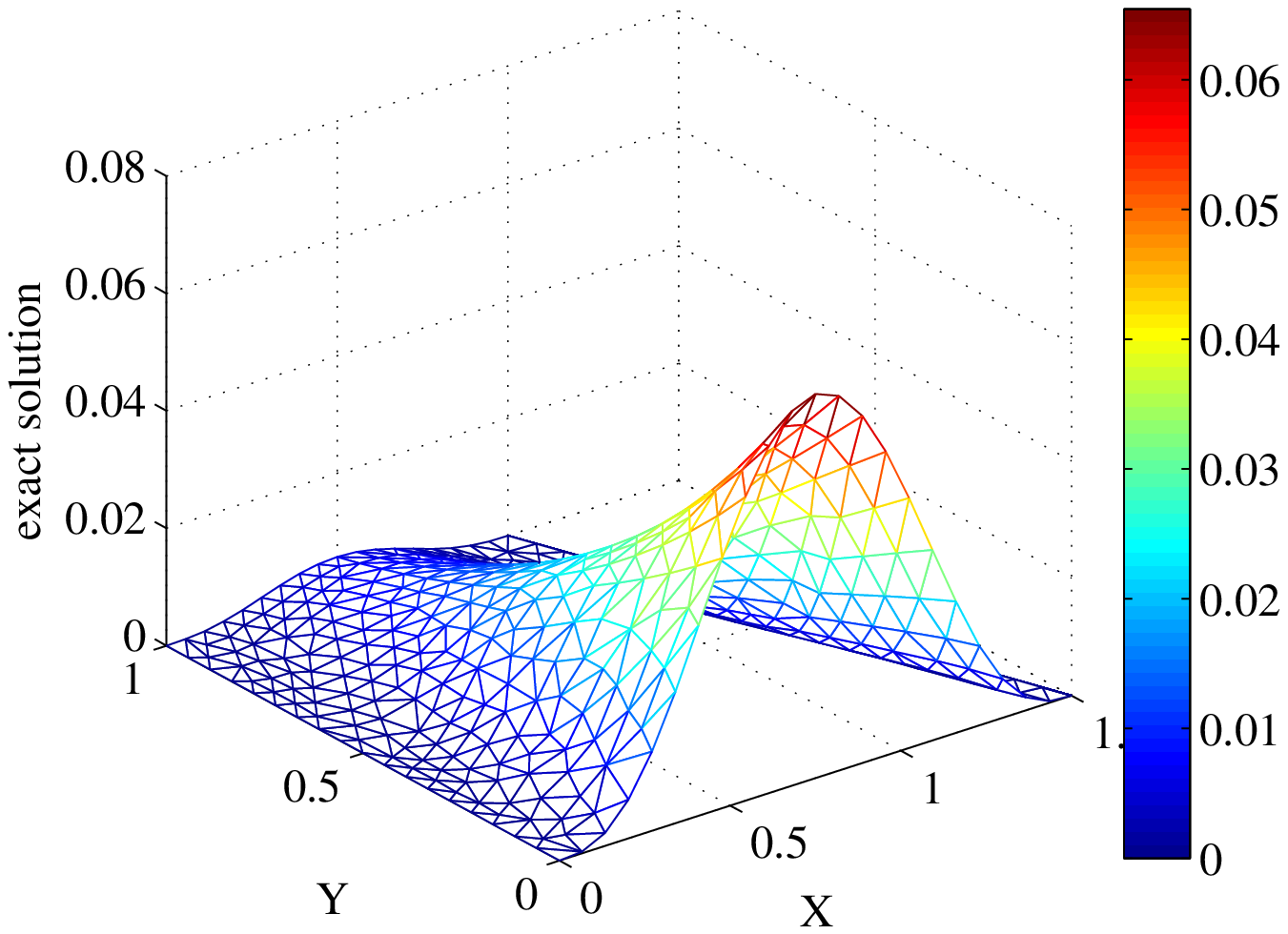}
\end{minipage}
\begin{minipage}[t]{0.49\linewidth}
\includegraphics[width=2.45in]{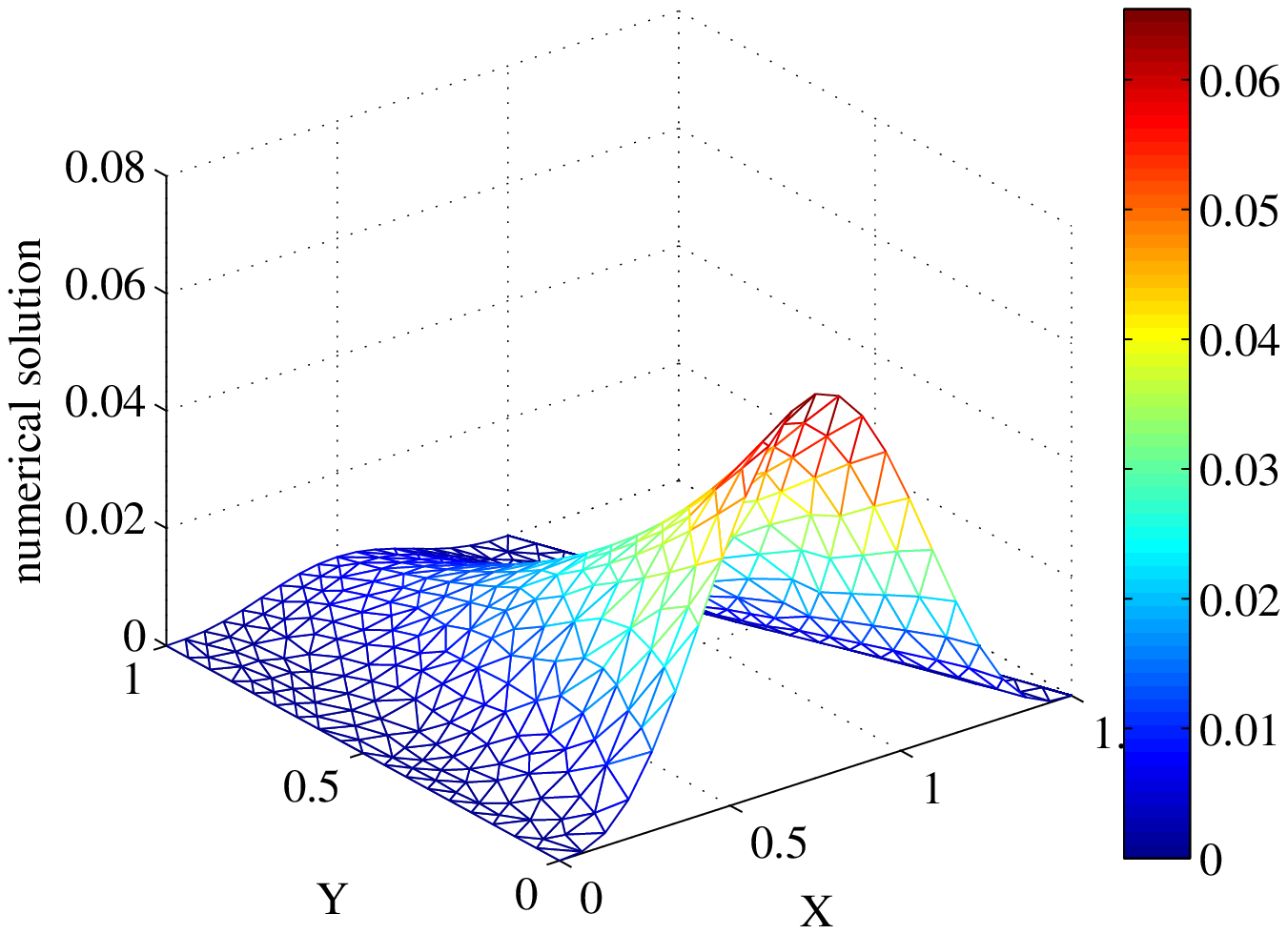}
\end{minipage}
\caption{The exact and numerical solutions created by IM-DQ method for Example 5.4.}\label{fig6}
\end{figure}

\noindent
\textbf{Example 5.4.} Consider the 2D fractional diffusion equation
\begin{align*}
\frac{\partial u(x,y,t)}{\partial t}-x^{\alpha_1} \mathcal{D}^{\alpha_1}_{\theta_1}u(x,y,t)-(1.5-x-0.5y)^{\alpha_2-3}\mathcal{D}^{\alpha_2}_{\theta_2}u(x,y,t)=f(x,y,t),
\end{align*}
on the trapezoidal domain as shown in Fig. \ref{fig5}, with $\theta_1=0$, $\theta_2=\pi$, the exact solution $u(x,y,t)=e^{-t}x^3( 0.5(3-y)-x )^3$, the source function
\begin{align*}
    f(x,y,t)= -e^{-t}x^3(0.5(3-y)-x)^3-e^{-t}(g^*_1(\alpha_1,x,3-y)+g^*_2(\alpha_2,x,y-3)),
\end{align*}
and the initial and boundary conditions taken from $u(x,y,t)$, where
\begin{align*}
    g^*_1(\alpha,x,y)&=\frac{0.75x^{3}y^3}{\Gamma(4-\alpha )}-\frac{18x^{4}y^2}{\Gamma(5-\alpha )}+\frac{180x^{5}y}{\Gamma(6-\alpha)}-\frac{720x^{6}}{\Gamma(7-\alpha)},\\
    g^*_2(\alpha,x,y)&=\frac{0.75(\alpha-2)(\alpha-1)\alpha y^3+18(\alpha-1)\alpha xy^2+180\alpha x^2y+720x^3}{\Gamma(7-\alpha)}.
\end{align*}
Let $\alpha_1=1.1$, $\alpha_2=1.3$, and $N=5000$. We run the algorithm with $c^*=0.75$ for Multiquadrics and $c^*=1.05$ for Inverse Multiquadrics
by using the scattered nodal points of total numbers 66, 171, 287, and 437, which are plotted in Fig. \ref{fig5}, respectively.
Table \ref{tab5} reports the errors of the numerical solutions of DQ methods with respect to the exact solution at $t=1$ at length.
It is evident that, under the chosen $\epsilon$'s, the proposed methods are stable and convergent on this trapezoidal domain.
Using the last nodal distribution plotted in Fig. \ref{fig5} and its corresponding $\epsilon$, a comparison between the exact and numerical solutions at
$t=1$ created by IM-DQ method is presented in Fig. \ref{fig6}, where no obvious difference can be observed from these two graphs.\\

\begin{table*}
\centering
\caption{The numerical results  at $t=1$ with $N=5000$, $\alpha_1=1.1$, and $\alpha_2=1.3$ for Example 5.4.}\label{tab5}
 \begin{tabular}{lclclc}
\toprule
\multicolumn{1}{l}{\multirow{2}{0.6cm}{$M$}} &\multicolumn{2}{l}{ MQ-DQ method} &\multicolumn{2}{l}{IM-DQ method} \\
\cline{2-5}& $e_2(\tau,M)$  &$e_\infty(\tau,M)$ &$e_2(\tau,M)$  &$e_\infty(\tau,M)$ \\
\midrule 65     &  2.3564e-04  &6.9019e-04   &3.1638e-04  &9.7823e-04 \\
         170    &  1.8822e-04  &6.3616e-04   &2.2830e-04  &7.6927e-04 \\
         286    &  1.0976e-04  &4.3639e-04   &1.3165e-04  &5.2928e-04 \\
         436    &  6.9543e-05  &2.7613e-04   &8.4616e-05  &3.5031e-04 \\
\bottomrule
\end{tabular} 
\end{table*}

\begin{figure}[!htb]
\begin{minipage}[t]{0.49\linewidth}
\includegraphics[width=2.75in]{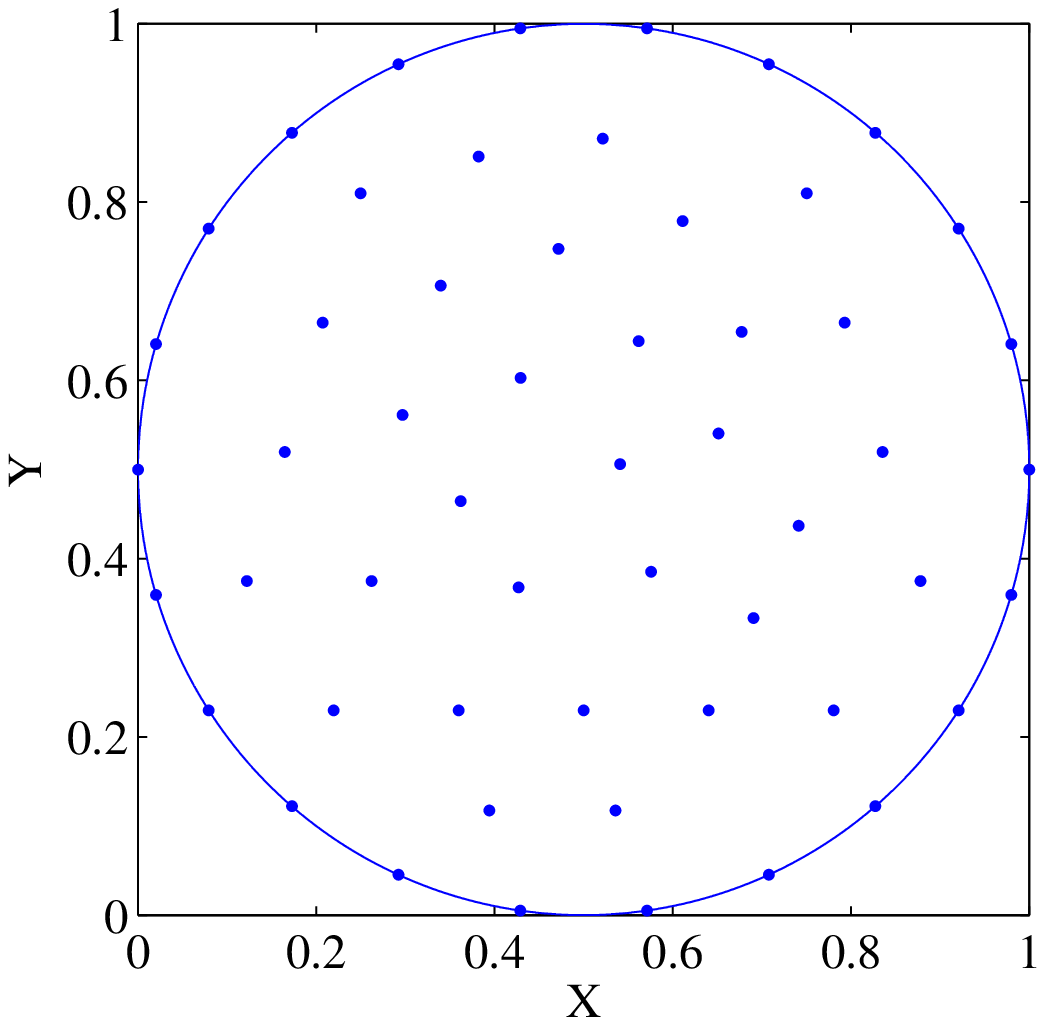}
\end{minipage}
\begin{minipage}[t]{0.49\linewidth}
\includegraphics[width=2.75in]{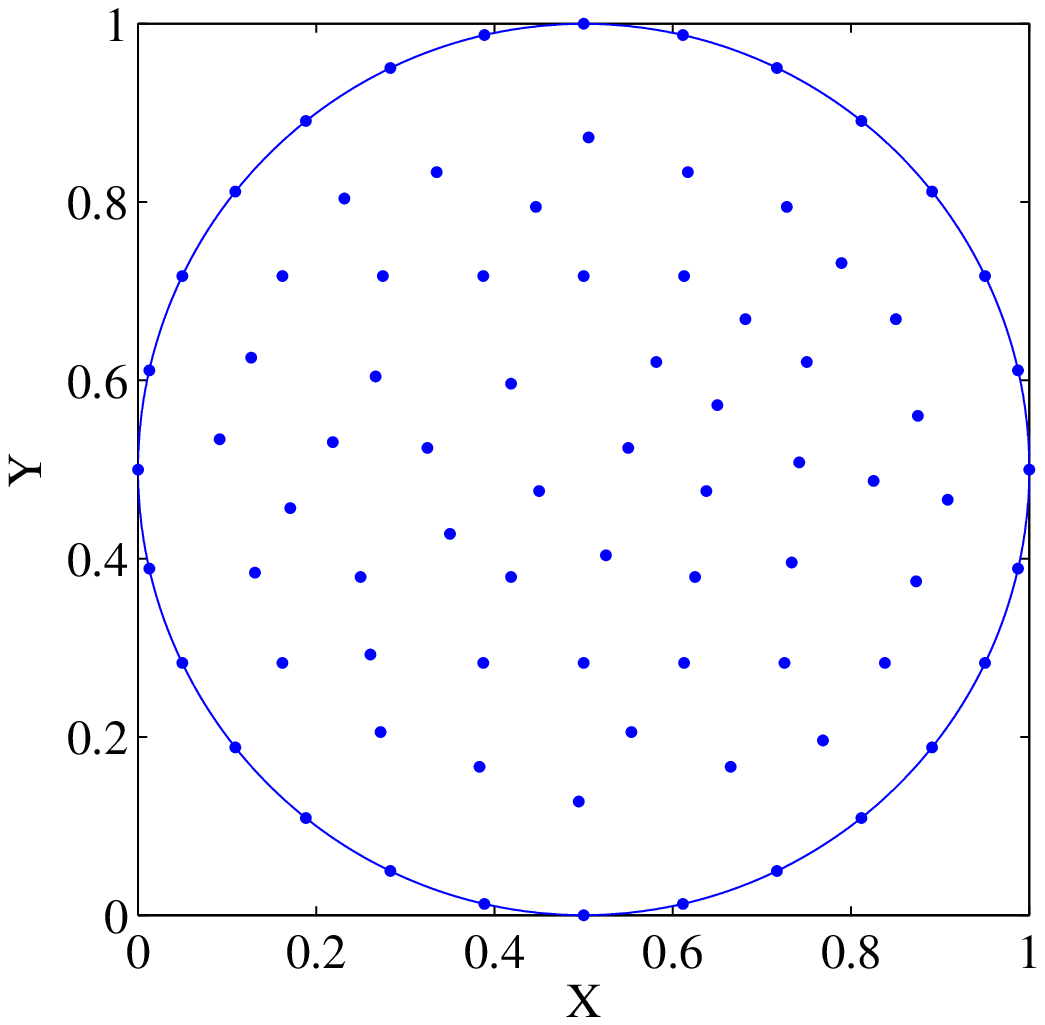}
\end{minipage}\\
\begin{minipage}[t]{0.49\linewidth}
\includegraphics[width=2.75in]{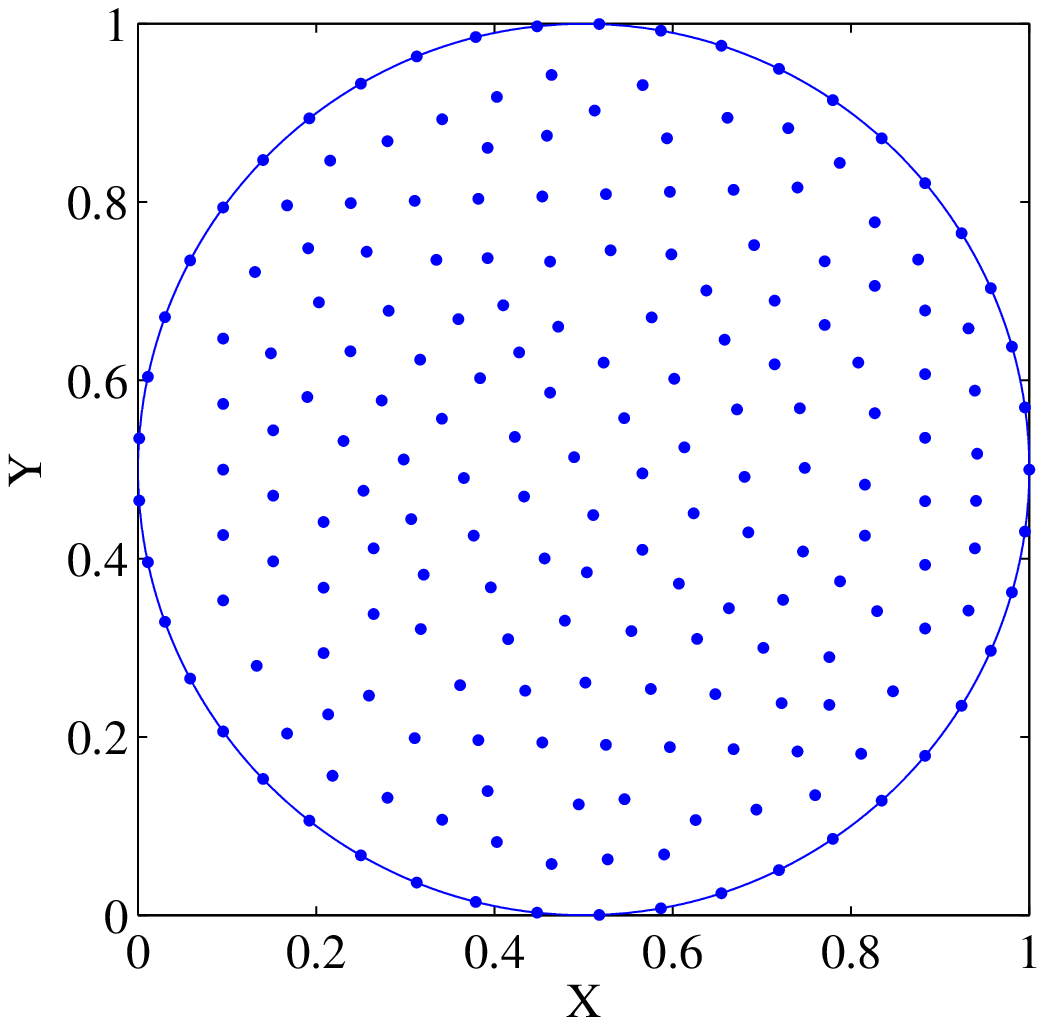}
\end{minipage}
\begin{minipage}[t]{0.49\linewidth}
\includegraphics[width=2.75in]{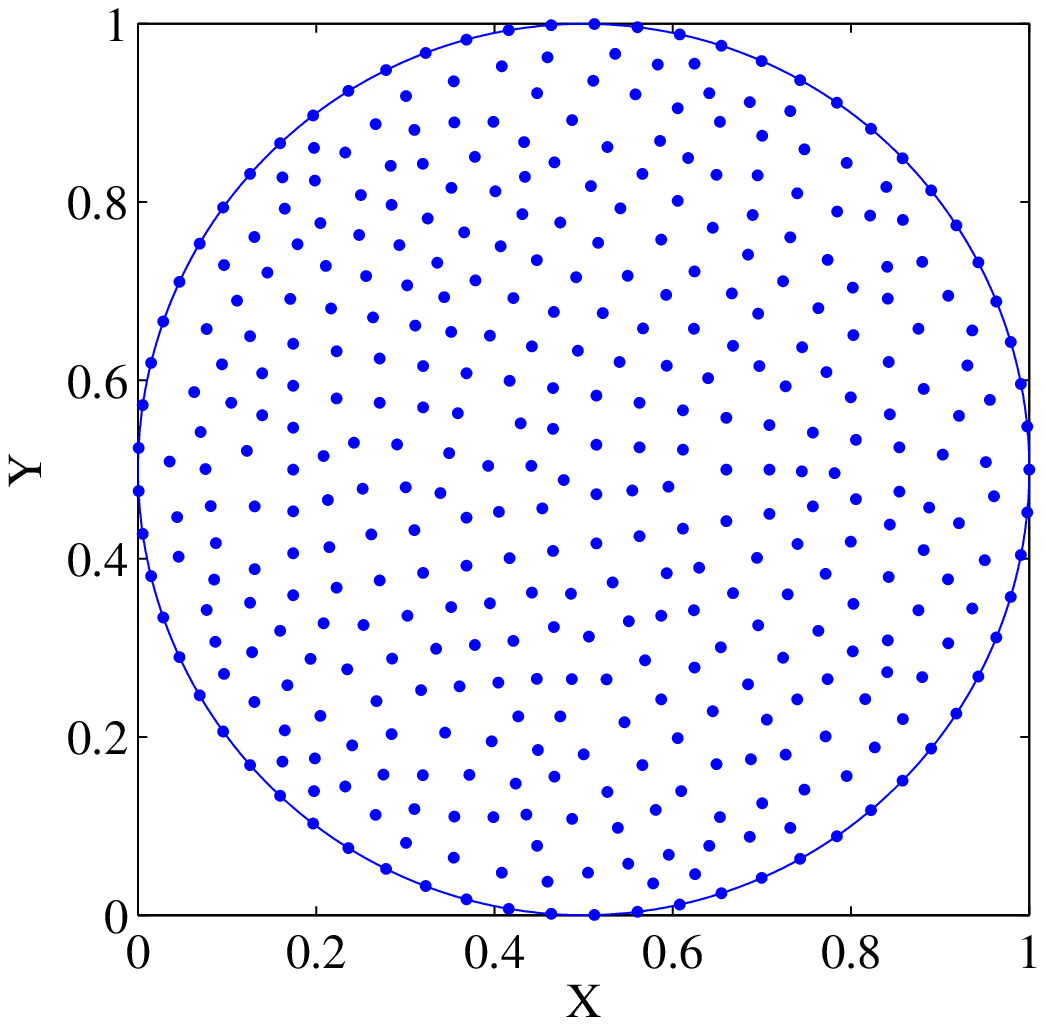}
\end{minipage}
\caption{The nodal distributions with nodal numbers $54$, $80$, $201$, and $402$ for Example 5.5.}\label{fig7}
\end{figure}

\begin{figure}[!htb]
\begin{minipage}[t]{0.49\linewidth}
\includegraphics[width=2.45in]{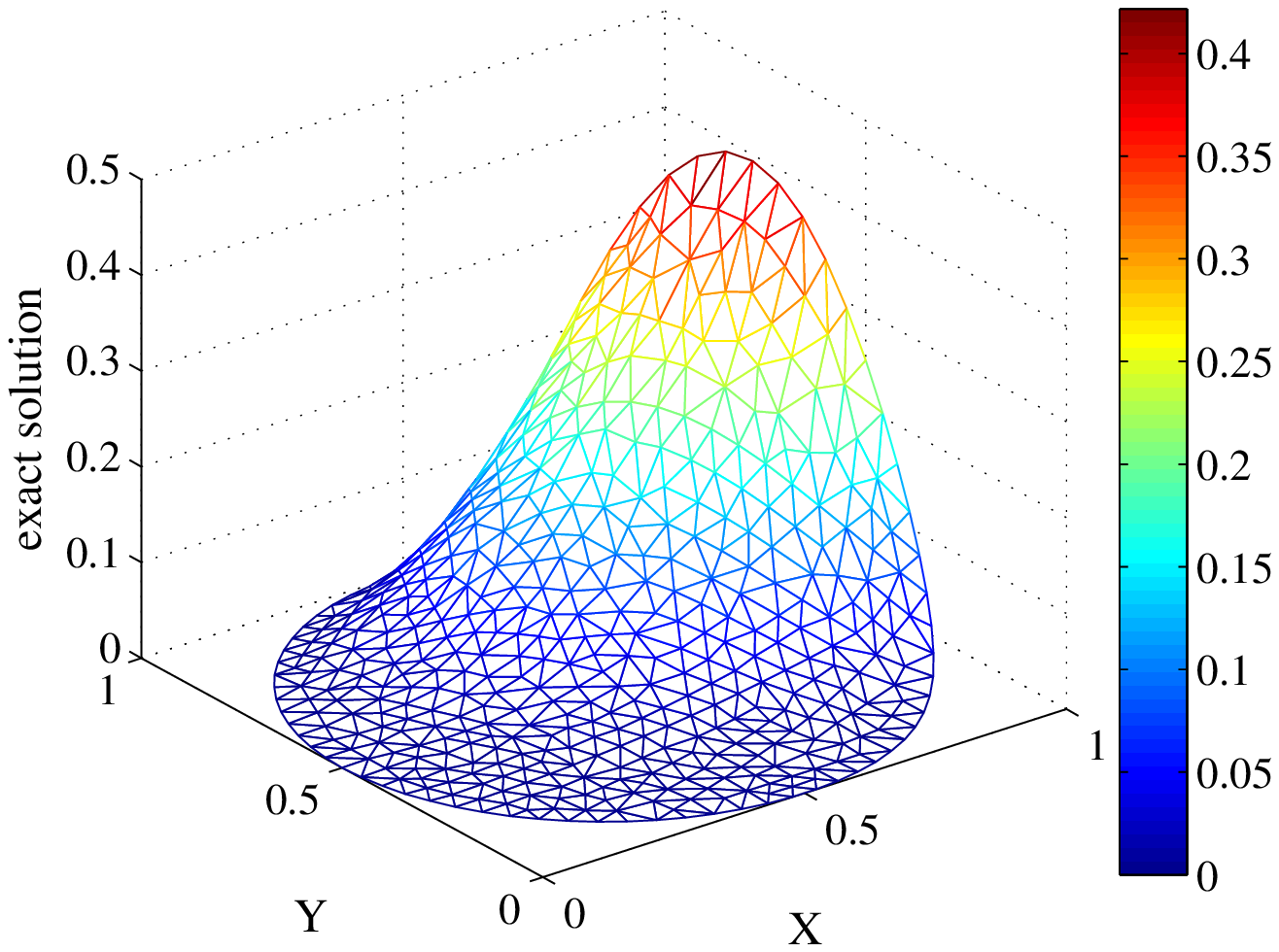}
\end{minipage}
\begin{minipage}[t]{0.49\linewidth}
\includegraphics[width=2.45in]{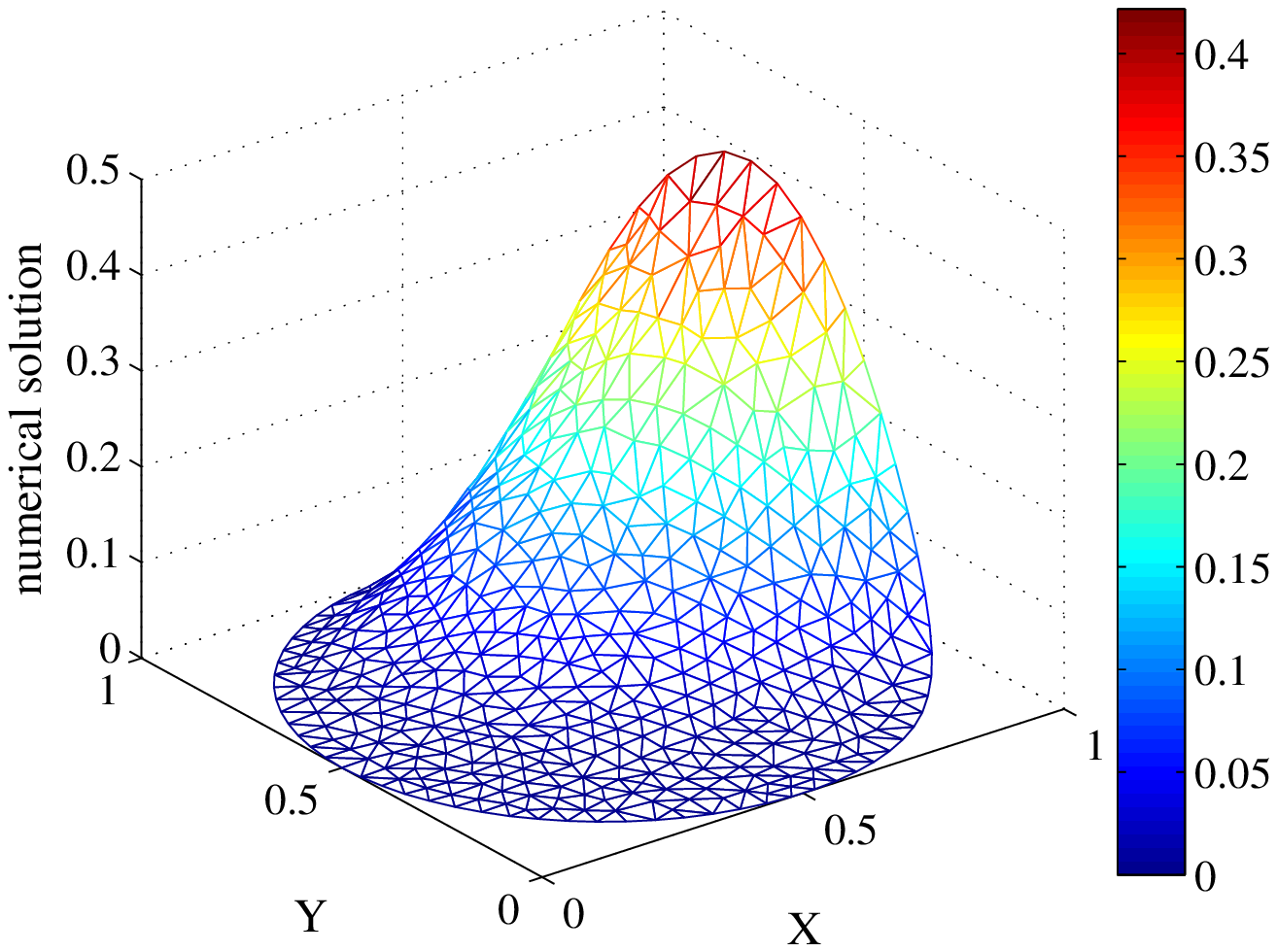}
\end{minipage}
\caption{The exact and numerical solutions created by GA-DQ method for Example 5.5.}\label{fig8}
\end{figure}

\noindent
\textbf{Example 5.5.} In this test, consider the 2D fractional diffusion equation
\begin{align*}
\frac{\partial u(x,y,t)}{\partial t}-\frac{y^\alpha}{2} \mathcal{D}^{\alpha}_{\theta}u(x,y,t)=f(x,y,t),
\end{align*}
on $\Omega=\{(x,y):(x-0.5)^2+(y-0.5)^2\leqslant0.25\}$, with $\theta=0$, $\alpha=1.9$,
the exact solution $u(x,y,t)=t^2(x-0.5+\sqrt{0.25-(y-0.5)^2})^2y^2$, the source function
\begin{align*}
    f(x,y,t) &= 2t(x-0.5+\sqrt{0.25-(y-0.5)^2})^2y^2 \\
       & \quad -\frac{t^2(x-0.5+\sqrt{0.25-(y-0.5)^2})^{2-\alpha}y^{2+\alpha}}{\Gamma(3-\alpha)},
\end{align*}
and the initial and boundary conditions taken from the exact solution.
The computational domain as well as the involved nodal distributions with nodal numbers 54, 80, 201, and 402
are displayed in Fig. \ref{fig7}, respectively. To see the actual performance of our methods, we compute the numerical results at $t=1$ with $N=5000$ and
document them side by side in Table \ref{tab6} by using $c^*=0.85$ for Inverse Multiquadrics and the data sets consisting of
5.4216, 5.9814, 7.5306, and 8.9554 for the shape parameter contained in Gaussians.
It is found that our methods show their capability to deal with the fractional
problem on circular domains even if the nodal distributions are irregular. Besides, under the given $\epsilon$'s, IM-DQ method
outperforms GA-DQ method in term of overall accuracy. Resetting $M=472$ together with $\epsilon=10.3221$,
we compare the exact and numerical solutions at $t=1$ created by GA-DQ method in Fig. \ref{fig8}. It is observed that
our method produces the approximation with sufficiently small errors, which can not be visually distinguished
from the exact one. This confirms the effectiveness of the proposed methods once more.

\begin{table*}[!htb]
\centering
\caption{The numerical results at $t=1$ with $N=5000$ and $\alpha=1.9$ for Example 5.5.}\label{tab6}
\begin{tabular}{lclclc}
\toprule
\multicolumn{1}{l}{\multirow{2}{0.6cm}{$M$}} &\multicolumn{2}{l}{ IM-DQ method} &\multicolumn{2}{l}{GA-DQ method} \\
\cline{2-5}& $e_2(\tau,M)$  &$e_\infty(\tau,M)$ &$e_2(\tau,M)$  &$e_\infty(\tau,M)$ \\
\midrule 53     &  4.4502e-03  &2.1437e-02   &1.2039e-02  &5.8457e-02 \\
         79     &  2.9459e-03  &1.3023e-02   &8.8637e-03  &3.8900e-02 \\
         200    &  7.3905e-04  &4.0762e-03   &1.8664e-03  &1.0135e-02 \\
         401    &  3.8098e-04  &2.3782e-03   &8.9460e-04  &6.1255e-03 \\
\bottomrule
\end{tabular}
\end{table*}

\begin{figure}[!htb]
\begin{minipage}[t]{0.49\linewidth}
\includegraphics[width=2.6in]{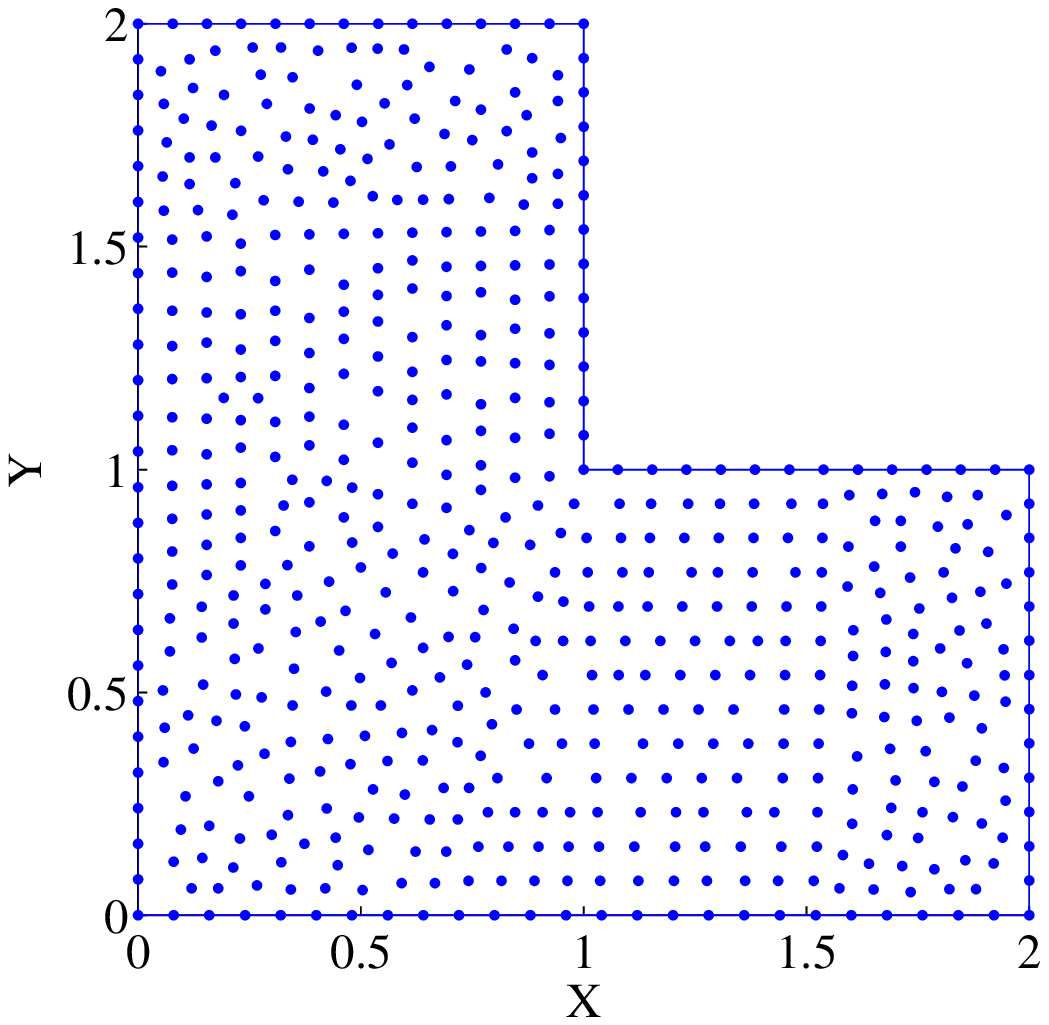}
\end{minipage}
\begin{minipage}[t]{0.49\linewidth}
\includegraphics[width=2.75in]{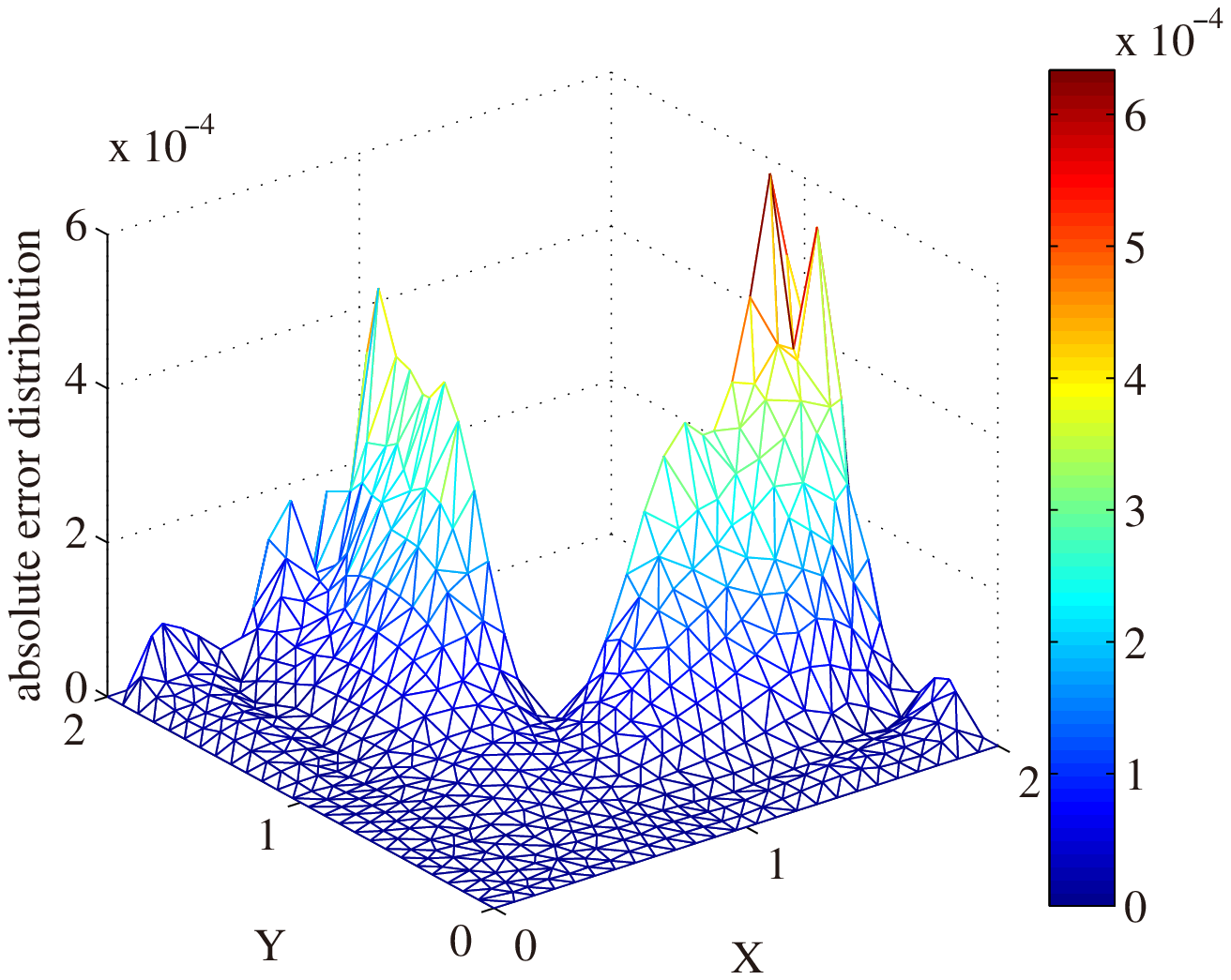}
\end{minipage}
\caption{The nodal distribution with $M=592$ and the absolute error distribution of IM-DQ method.}\label{fig9}
\end{figure}

\noindent
\textbf{Example 5.6.} In the last test, consider the 2D fractional diffusion equation
\begin{align*}
\frac{\partial u(x,y,t)}{\partial t}-x^\alpha y^\alpha\sum^3_{l=1}\mathcal{D}^{\alpha}_{\theta_l}u(x,y,t)=f(x,y,t),
\end{align*}
on the L-shaped domain as shown in Fig. \ref{fig9},  with $\theta_1=0$, $\theta_2=\pi/4$, and $\theta_3=\pi/2$,
subjected to the initial and boundary conditions taken from the exact solution $u(x,y,t)=t^3x^2y^2$. The source function is manufactured by
\begin{align*}
    f(x,y,t) &= 3t^2x^2y^2 -t^3x^\alpha y^\alpha\bigg(f^*(x,y,t)+2\frac{x^{2-\alpha}y^2+x^2y^{2-\alpha}}{\Gamma(3-\alpha)}\bigg),
\end{align*}
with $\alpha_1=\alpha$ in $f^*(x,y,t)$ declared above. Letting $\epsilon$ be 0.2128 in Multiquadrics, 0.3445 in Inverse Multiquadrics,
and 4.6880 in Gaussians, we examine the convergence of DQ methods at $t=0.5$ versus the variation of $\alpha$ with $N=1000$ on
the scattered nodal points of total number $593$. The numerical results are reported in Table \ref{tab7}.
The used irregular nodal distribution and the absolute error distribution of IM-DQ method for $\alpha=1.5$ are plotted in Fig. \ref{fig9}, respectively.
As one can see, the errors decay as $\alpha$ increases in addition to a special case of $\alpha=2.0$ for $e_\infty(\tau,M)$, and
MQ-DQ method outperforms IM-DQ, GA-DQ methods in term of overall accuracy under the given $\epsilon$'s.

\begin{table*}
\centering
\caption{The numerical results at $t=0.5$ with $N=2000$ and $M=592$ for Example 5.6}\label{tab7}
\begin{tabular}{lllllll}
\toprule
\multicolumn{1}{l}{\multirow{2}{0.5cm}{$\alpha$}}
&\multicolumn{2}{l}{ MQ-DQ method} &\multicolumn{2}{l}{ IM-DQ method} &\multicolumn{2}{l}{GA-DQ method} \\
\cline{2-7}& $e_2(\tau,M)$  &$e_\infty(\tau,M)$ &$e_2(\tau,M)$  &$e_\infty(\tau,M)$ &$e_2(\tau,M)$  &$e_\infty(t,M)$ \\
\midrule 1.2    &1.5847e-04  &5.3015e-04  &  1.4751e-04  &8.9850e-04 &2.9393e-04  &1.5306e-03 \\
         1.5    &1.0553e-04  &4.0805e-04  &  1.1669e-04  &6.3374e-04 &2.5013e-04  &1.1697e-03 \\
         1.8    &6.3716e-05  &3.6591e-04  &  8.9356e-05  &4.8952e-04 &1.9519e-04  &8.8947e-04 \\
         2.0    &5.2515e-05  &3.8395e-04  &  7.2907e-05  &4.4158e-04 &1.5855e-04  &7.1631e-04 \\
\bottomrule
\end{tabular}
\end{table*}

\section{Conclusion}
The multi-dimensional space-fractional diffusion equations are essential part of fractional
PDEs and have been one of principle concerns in mathematical physics,  but solving these types of equations
appears to be somewhat challenging due to the global correlation of fractional derivatives, especially on irregular domains.
In this research, by using Multiquadric, Inverse Multiquadric and Gaussian RBFs as trial functions, the DQ formulations
for fractional directional derivatives of Caputo type are presented. Then, three effective DQ methods are proposed for
the space-fractional diffusion equations on 2D irregular domains via discretizing the resultant ODEs by employing
the Crank-Nicolson scheme. These methods enjoy some appealing advantages such as low occupancy cost, truly \emph{mesh-free},
and the adaptability to arbitrary domains, which make they compare favorably to some of traditional methods as FDM.
The computational accuracy relies on the nodal number, shape parameter, the complexity of fractional
models, and the other potential factors. The codes are utilized to treat the fractional problems
on square, trapezoidal, circular, and L-shaped domains, respectively, and the outcomes have demonstrated that
our methods are capable of capturing the exact solutions under the condition that the shape parameters are well prepared.
Due to the characteristics of RBFs, these DQ methods are insensitive to dimensional
change, which means that they can be generalized to three-dimensions without
a significant increase in computational burden. This merit evidently offers us a possibility
to simulate the other space-fractional problems in high-dimensions with complex boundaries
arising in the fields of anomalous dispersion, fluid dynamics, conservative systems, genetic propagation, and viscoelastic materials.

\begin{acknowledgements}
The authors would like to thank the anonymous referees for their valuable comments and suggestions.
This research was supported by National Natural Science Foundations of China (Nos. 11471262 and 11501450).
\end{acknowledgements}

\bibliographystyle{spmpsci}      
\bibliography{mybib}   

\end{document}